\documentclass{article}[12pt]
\usepackage{color,times,amsmath,amsfonts,latexsym,epsfig,epsf,colordvi}
\usepackage{url,hyperref}
\usepackage[english]{babel}

\newcommand{\CC}{\mathbb {C}}
\newcommand{\RR}{\mathbb {R}}

\newcommand{\al}{\alpha }

\newcommand{\et}{\eta }

\newcommand{\bp}{\begin{pmat}}
\newcommand{\ep}{\end{pmat}}
\def\tcr{\textcolor[rgb]{1,0,0}}
\def\tcb{\textcolor[rgb]{0,0,1}}

\def\circledzero{{\begin{picture}(5.6,3.5) %
   \put(2.2,1.3){\circle{3.8}} %
   \put(1.3,.2){\normalsize 0} %
   \end{picture}}}   
\def\circledone{{\begin{picture}(5.6,3.5) %    <<<<====
   \put(2.2,1.3){\circle{3.8}} %
   \put(1.3,.2){\normalsize 1} %
   \end{picture}}}   
 \def\circledtwo{{\begin{picture}(5.6,3.5) %
   \put(2.2,1.3){\circle{3.8}} %
   \put(1.3,.2){\normalsize 2} %
   \end{picture}}}   
 \def\circledthree{{\begin{picture}(5.6,3.5) %
   \put(2.2,1.3){\circle{3.8}} %
   \put(1.3,.2){\normalsize 3} %
   \end{picture}}}      
   \def\circledfour{{\begin{picture}(5.6,3.5) %
   \put(2.2,1.3){\circle{3.8}} %
   \put(1.3,.2){\normalsize 4} %
   \end{picture}}}   
\def\circledfive{{\begin{picture}(5.6,3.5) %
   \put(2.2,1.3){\circle{3.8}} %
   \put(1.3,.2){\normalsize 5} %
   \end{picture}}}   
 \def\circledsix{{\begin{picture}(5.6,3.5) %  << == 
   \put(2.2,1.3){\circle{3.8}} %   was 1.5, 1.3
   \put(1.3,.2){\normalsize 6} %   was 0 , .3
   \end{picture}}}   
 \def\circledseven{{\begin{picture}(5.6,3.5) %
   \put(2.2,1.3){\circle{3.8}} %
   \put(1.35,.2){\normalsize 7} %
   \end{picture}}}
\def\circleddots{{\begin{picture}(5.6,3.5) %
   \put(2.2,1.3){\circle{3.8}} %
   \put(1.35,.25){\normalsize $..$} %
   \end{picture}}}
 \def\circledstar{{\begin{picture}(5.6,3.5) %
   \put(2.2,1.3){\circle{3.8}} %
   \put(1.35,.34){\normalsize $*$} %
   \end{picture}}}

\author{  Frank Uhlig \thanks{Department of Mathematics and Statistics, Auburn 
University, Auburn, AL 36849-5310 \ (uhligfd@auburn.edu)}}

 \title{\vspace*{-12mm} Adapted AZNN Methods for Time-Varying and Static Matrix Problems\\}
\begin{document}
\date{~}
\thispagestyle{empty}
\maketitle

\thispagestyle{empty}

\vspace*{-16mm}
\begin{center} { \bf Abstract  } \\[1mm]
\begin{minipage}{150mm}
We present  adapted Zhang Neural Networks (AZNN) in which the parameter settings for the exponential decay constant $\eta$ and the length of the start-up phase of basic ZNN are adapted to the problem at hand. 
Specifically, we study experiments with AZNN for time-varying square matrix factorizations as a product of time-varying symmetric matrices and for the time-varying matrix square roots problem. Differing from generally used small $\eta$ values and minimal start-up length phases in ZNN, we adapt the basic ZNN method  to work with large or even gigantic $\eta$ settings and arbitrary length start-ups using Euler's low accuracy finite difference formula. These adaptations improve the speed of AZNN's convergence and lower its solution error bounds for our chosen problems significantly to near machine constant or even lower levels.\\
 Parameter-varying AZNN also allows us to find full rank symmetrizers of static matrices reliably, for example for the Kahan and Frank matrices  and for matrices with highly ill-conditioned eigenvalues and complicated Jordan structures of dimensions from $n = 2$ on up. This helps  in cases where full rank static matrix symmetrizers have never been successfully computed before.
\end{minipage}\\[-1mm]
\end{center}  
\thispagestyle{empty}



\pagestyle{myheadings}
\thispagestyle{plain}
\markboth{Frank Uhlig}{Adapted AZNN }

\section{Introduction to AZNN Methods for  Real or Complex  Matrix Flows $A(t)$}

\vspace*{-1mm}

We consider time-varying matrix problems that we want to solve more speedily and accurately by using Zhang Neural Networks (ZNN) adaptively and judiciously, and varying the length of the initial start-up phase and also adjusting the exponential decay constant $\eta$ in  AZNN  to the three phases of ZNN, namely its start-up, middle, and final stages. Zeroing Neural Networks such as as Zhang Neural Networks are  computational workhorses in many branches of engineering. For functional inputs, discretized ZNN helps us to understand  the feasibility of a given parametric matrix model's function equation over time or learn of its unfeasibility. If feasible, the chosen model can be implemented successfully  in ZNN for sensor clocked discretized input data and used on-chip in robots, other machinery, and in other engineering processes.\\[1mm]
 ZNN goes back to Yunong Zhang's Ph.D. work of 2000/2001 in Hong Kong \cite{ZW01}. ZNN's engineering applications have been described in almost 500  papers  as well as in a handful of  books. A theoretical numerical analysis of ZNN methods has never been attempted, see \cite{FUsurveyZNN} e.g..\\[1mm]
Starting from a time-varying matrix/vector problem\\
 \hspace*{12mm}\circledzero \ \ $f(x(t),A(t), ...) = g(t, ...)$, \\
 discretized ZNN is designed to compute an approximate solution iteratively in time  from a few start-up values. It does so  predictively over time and with high accuracy. ZNN's success hinges on Zhang's idea of stipulating exponential decay for the associated error function\\
\hspace*{12mm}\circledone \ \ $E(t) = f(x(t),A(t),...) - g(t,...)$,\\
 namely \\
 \hspace*{12mm}\circledtwo \ \ $\dot E(t) = - \eta E(t)$ for $\et > 0$.\\
\hspace*{12mm}\circledthree \ \ Assuming  a constant sampling gap $\tau = t_{k+1} - t_k = const$ throughout the ZNN process,  the derivative of the unknown $\dot x(t_k)$ in the error differential  equation  \circledtwo is isolated algebraically if possible. This typically leads to an expression for  $\dot x(t_k)$ as the solution of a system of linear equations 
$$
 \dot x(t_k) = P(t_k) \backslash q(t_k) .
 $$
  \hspace*{12mm}\circledfour \ \  After summing  several  Taylor expansions  around $x(t_k)$ that contain future $x(t_{k+1})$  and current or past solution $x(t_\ell)$ data for $\ell \leq k$ and the derivative $\dot x(t_k)$ of the unknown, we use a look-ahead and convergent finite difference formula of type {\tt j\_s}, see \cite{FUfindiff19} for a list thereof, and solve for $\dot x(t)$. This leads to an expression for the derivative $ \dot x(t_k)$ such as 
  $$ 
  \dot x_k = \dfrac{8x_{k+1}+x_k -6x_{k-1} - 5 x_{k-2} + 2 x_{k-3}}{18\tau}
  $$
 when for example using a  five instance finite difference  formula of type {\tt j\_s} = {\tt 2\_3} with global truncation error order $O(\tau^3)$. \\
 \hspace*{12mm}\circledfive \ \  Then  the two expressions of the derivative  $\dot x(t_k)$ in \circledthree and \circledfour of the unknown $x(t)$ are equated.\\
   \hspace*{12mm}\circledsix \ \  The resulting solution-derivative free equation of \circledfive is finally solved  for $$x(t_{k+1}) =  \dfrac{9}{4} \tau (P(t_k) \backslash q(t_k)) - \dfrac{1}{8} x_k + \dfrac{3}{4} x_{k-1} + \dfrac{5}{8} x_{k-2} - \dfrac{1}{4} x_{k-3} .
   $$
    This predicts the problem's solution at the next time step $t_{k+1}$. \\
 \hspace*{12mm}\circledseven \ \  Finally ZNN   iterates the process.\\[1mm]
 In step \circledsix$\!\!$, the discretized ZNN method solves one linear equation and evaluates one short recursion of past solution data, and thereby it solves the underlying time-varying matrix/vector problem cheaply and predictively in time.\\
 
   The above seven \circleddots steps of ZNN have recently been derived and explained in detail inside \cite{FUsurveyZNN}.\\[1mm]
ZNN methods may seem connected to analytic continuation methods and  ODE solvers when they start with the error differential equation \circledtwo$\!\!$, but when interpreting  their operations, results and watching their phenomena, ZNN methods appear to follow  quite different numerical principles and seem to belong to a different area of numerical matrix analysis. Further details, open questions,   and an early assessment of ZNN methods are given \cite{FUsurveyZNN}.\\[1mm]
This paper investigates both the role of the exponential decay constant $\eta$ of ZNN and the actual length of ZNN’s start-up phase. Our adapted ZNN method (called AZNN) adopts different strategies to improve  the convergence behavior of ordinary ZNN. Adaptive $\eta$ and other settings help AZNN  compute the time-varying problem's predictive solution more accurately and much sooner than ZNN. This is specially useful for time-varying matrix problems for which no globally reliable static matrix algorithms exist such as for the factorization of  real or complex square  matrices into the product of two symmetric ones \cite{DU16}. Note  that ZNN methods are very tolerant of random, partially 'little-sense' start-up values. But such start-up value settings often  lead to divergence of the subsequent iterative process or they let ZNN  needlessly wander through 'solutions' with very low accuracy  for a long time.  A better way to obtain reasonable start-up values is to use Euler's look-ahead method of type {\tt j\_s = 1\_2} with local truncation error order $O(\tau^2)$\\[-2mm]
$$  
x(t_{k+1}) = x(t_k) + \tau \cdot \dot x(t_k) + O(\tau^2) \ \   \text{ \ that is based on the secant rule} \ \ \ \dot x(t_k) = \dfrac{x(t_{k+1} ) - x(t_k)}{\tau} + O(\tau) ,
$$

\vspace*{-1mm} 
 and the problem formulation in tandem from an initial  random first guess for $x({t_0})$ and with an adapted $\eta$ setting. In contrast, the basic non-adapted ZNN method  uses the identical value of $\eta$ as the exponential decay constant  in its start-up phase and in its main ZNN iterations \circledfive through \circledseven$\!\!$.  Likewise, frugality concerns for basic ZNN suggest to compute  only the minimally needed number of start-up data before switching into iterative  ZNN mode. \\[1mm]
In AZNN we use different $\eta$ settings in the start-up point generation phase and  in the iterations phase and we test  longer and shorter sequences of start-up value runs to find the best option for  AZNN and the given matrix problem.  The Euler start-up run should last until its errors reach the Euler method's truncation error plateau of around $10^{-5}$ or $10^{-6}$, or a little bit beyond. Then  we   switch to the AZNN iterations phase where we use  a different $\eta$ value. These two adaptations achieve the chosen finite difference formula's  truncation error order plateau faster  than  standard ZNN does. The AZNN method helps us to solve our experimental example runs for finding time-varying matrix square roots and time-varying  symmetric matrix factorizations in a better way.\\[2mm]
Section 2 deals pragmatically with techniques to alter or adapt the basic ZNN method to two  specific parametric matrix equations.  Section 2 (A)  deals with the time-varying matrix square root problem $A(t) = X(t) \cdot X(t)$ that generally has many solutions unless its Jordan structure is badly formed such as in the 2 by 2 flow matrix $A(t) = \bp 0&t\\0&0\ep$ for which there exists no matrix square root at all for any $t \neq 0$, see \cite{EU92} or \cite{CL74}. In Section 2 (B) we consider the time-varying matrix symmetric factorization problem $A(t) = S(T) \cdot V(t)$ for two symmetric matrix flows $S(t) = S^T(t) $ and $V(t) = V^T(t) \in \CC_{n,n}$ with $S(t)$ nonsingular for all $t$. In theory, the symmetrizer problem for static entry matrices $A\in \CC_{n,n}$ can always be solved by solving the linear system in the entries of the symmetric matrix $S$ with $S\cdot A = (S\cdot A)^T = A^T \cdot S^T$. But this intuitive approach generally leads to very ill conditioned linear systems. The study of symmetric matrix factorizations  goes back more than a century to Frobenius  \cite[p. 421]{F1910}. Frobenius  showed that the symmetric  factor $S$ can always be chosen nonsingular and that the symmetrizers of $A$ form a linear subspace  of dimension $n$   for static $A_{n,n} \in \CC_{n,n}$. \\[1mm]
Static matrix symmetrizing algorithms that use linear equation solvers, iterative methods, or eigenvalue and singular value based methods were studied in depth in \cite{DU16}. All of these methods  were shown to be deficient in \cite{DU16} when a given matrix $A$ had degenerate  eigen or  Jordan structures. The  analysis in \cite{DU16} showed that eigen ill-conditioning  prevented these static methods generally from  finding invertible symmetrizers $S= S^T$ with $S\cdot A = V$ and $V = V^T$.  Until now no symmetric factorizations $A = S^{-1} \cdot V$ could be computed reliably for  such $A$. Section 3 of this paper finally shows how the adapted AZNN method can succeed via a simple parameterization in a  matrix flow related to such $A$ where static matrix theory and static numerical analysis have failed us before.\\[1mm]
{\bf Note :} The matrix symmetrizer problem deals exclusively with real or complex symmetric matrices $S = S^T$ and $V = V^T$ whose entries in position $i,k$ are identical to entries in position $k,i$ for all row and column indices $1 \leq i,k \leq n$. Matrices with complex conjugate entries in positions $i,k$ and $k,i$ are hermitean and not the subject of symmetric matrix factorizations, nor of this paper.\\[-6mm]

\section {Developing AZNN methods for time-varying square roots $X(t)$ and symmetric factorizations $S(t) \cdot V(t)$ of time-varying matrix flows $A(t)$} 

{\large \bf (A) \ \ \   Time-varying square roots $X(t)$ for time-varying matrix flows $A(t)$}\\[2mm]
\enlargethispage{10mm}
Fixed entry matrix $A_{n.n}$ and their square roots $X$ with $A = X \cdot X$ have been studied for at least half a century, both numerically and theoretically, see \cite[p. 466, 469, 506]{EU92} and \cite{CL74} e.g. for the existence or non-existence of static matrix square roots. We recommend the {\em Wikipedia} entry on {\em Square root of a matrix} and a look at numerical and applied methods. Use Matlab's {\tt sqrtm.m} file for computing $X$.\\[1mm]
Computing time-varying matrix square root flows $X(t)$  for time-varying matrix flows $A(t)$ requires more effort.\\ 
 Here we first look at the basic ZNN method for a time-varying matrix square root problem and recreate \cite[Figure 1]{FUsurveyZNN} for further study. In Figure 1 below, the two graphs show the error function value at the start-up to full {\tt j\_s = 4\_5} switching point displayed by a red square \tcr{$\square$} and the 1000th iterate's error function value displayed by a blue circle around a star \tcb{\circledstar}$\!\!$; on the left hand side with $\eta = 1.35$ as in \cite[Figure 1]{FUsurveyZNN}, and on the right side with $\eta = 2.6$ but for the higher truncation order finite difference formula of type {\tt j\_s = 4\_5}. Both switching points between the start-up and iteration phases of ZNN  occur at around 20 seconds, but with very different error function values, namely at around $ 4 \cdot 10^{-12}$ when $\eta = 1.35$ and at  around $4 \cdot 10^{-3}$ when $\eta = 2.6$. With $\eta = 1.35$ in the left hand graph of Figure 1, the transition from the start-up phase that creates $j+s = 4 + 5 = 9$ data points via Euler's method is very smooth until the Euler method's error function hovers below $ 10^{-11}$. For $\eta = 2.6$ the transition from Euler to basic $\eta = 2.6$ ZNN on the right of Figure 1 is also rather smooth for 20 + iteration steps, but thereafter  many steps are wasted with error function value going up and down until the {\tt 4\_5} based ZNN method stabilizes and takes over smoothly to converge well error-wise to below $10^{-11}$ after around 70 seconds from start as Figure 1 right and the enlarged view in Figure 2  clearly shows.\\[1mm]
% use ... 
%  /Box/local/matlabin/AZNN/tvMatrSquareRootwEulerStartva.m  in fig 1
\hspace*{-1.5mm}\includegraphics[height=74mm]{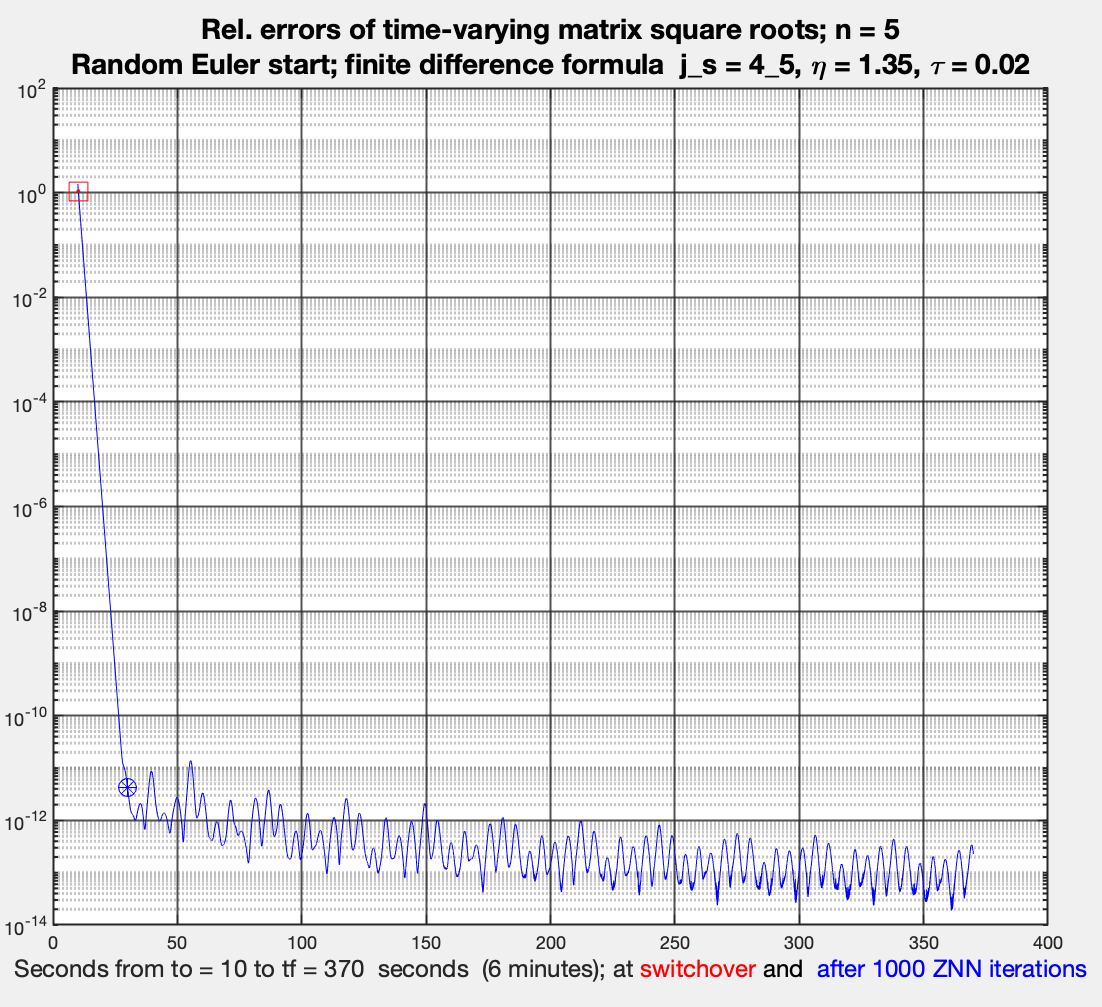} \vspace*{-2mm}\includegraphics[height=74mm]{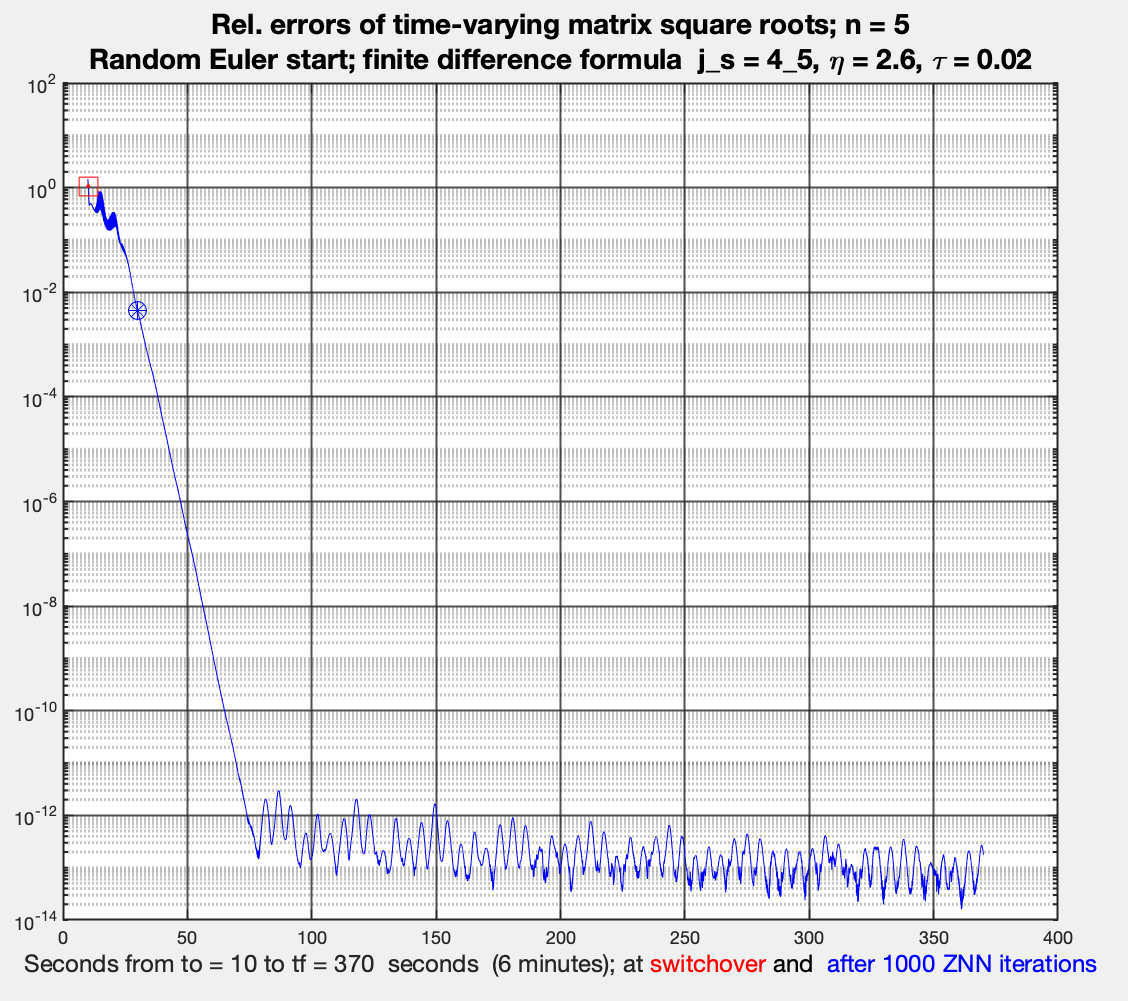}
\\[1mm]
\hspace*{73mm} Figure  1

If we increase $\eta$ further above 2.6, the error curve will move more or less horizontally around the $10^{-1}$ error level and  basic ZNN becomes useless. The error function magnitudes may even diverge to $10^{150}$ and higher values over our 6 minutes run.\\ % keep old nongrided image below
\vspace*{-6mm}
\begin{center}
\includegraphics[width=101mm]{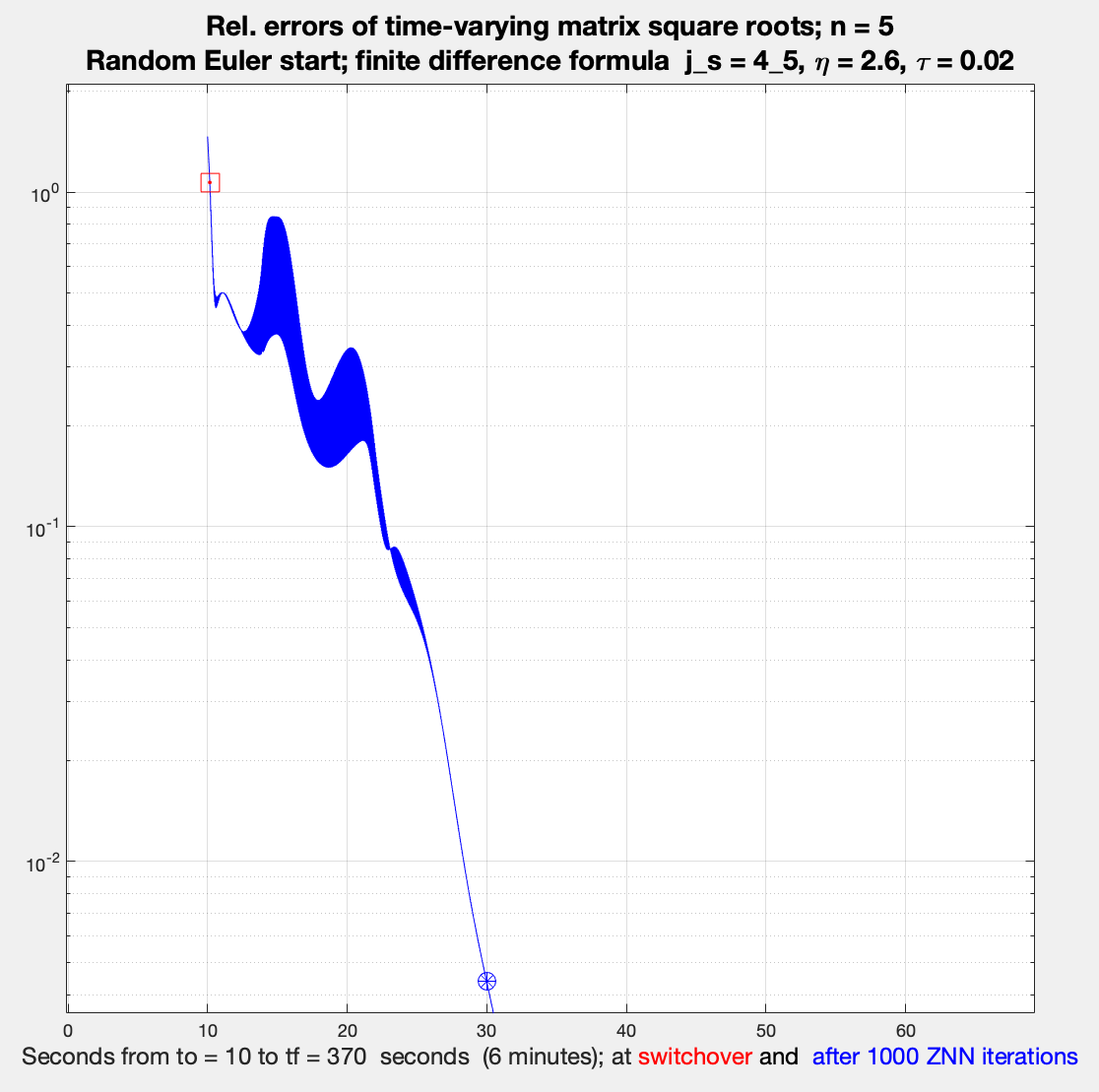}\\[-1mm]
Figure  2 
\end{center}

\vspace*{-2mm}
The aim of creating an adapted ZNN method is to try and extend the short region of rapidly  decreasing errors after the switch from Euler to full $j\_s$ ZNN steps depicted in Figure 2 and at the same time try to lower the terminal error function values to the theoretical $O(\tau^{j+s+1})$ local truncation error level for $j\_s$ based ZNN iterations.\\[1mm]
To understand how this can be achieved, we shall explain the seven steps of ZNN briefly now.\\[1mm]
\enlargethispage{5mm}
These are the seven steps in basic discretized ZNN for the time-varying matrix square roots problem for a time-varying matrix flow $A(t)_{n,n}$ with $ t \in [to,tf] \subset \RR$ and a constant sampling gap $\tau = t_{\ell + 1} -  t_\ell$ for all $\ell = 1,2,... $ .\\
\circledzero \ Model equation  \ $ A(t) = X(t) \cdot X(t)$ for a given matrix flow $A(t)_{n,n}$ and its unknown square roots $X(t)_{n,n}$.\\
\circledone \ Error equation \ $ E(t) = A(t) - X(t) \cdot X(t)$.\\
\circledtwo \ Differentiated  error equation \ $\dot E(t) = - \eta E(t)$.\\
\hspace*{5mm} \ Error DE reordered \ $\dot A(t) - \dot X(t)\cdot X(t)  - X(t) \cdot \dot X(t) = - \eta (A(t) - X(t) \cdot X(t))$; \ and by leaving off the \\
\hspace*{6.6mm}time $t$ dependence\\
\circledthree \hspace*{40mm}$\dot X X + X \dot X = \dot A + \et~ (A - X X)$.\\ 
\hspace*{5mm} \ Rephrased using  Kronecker products \ $A \otimes B$\  and the matrix column vector notation \ $X(:)$ \\ \circledthree$\!\!\!_K$ \hspace*{15mm}$(X^T \otimes I_n + I_n \otimes X)_{n^2,n^2}\cdot  \dot X(:)_{n^2,1} \in \CC^{n^2} =$ \\
\hspace*{40mm}$= \dot A(:)_{n^2,1} + \et\ A(:)_{n^2,1} - \et\  (X^T \otimes I_n)_{n^2,n^2} \cdot X(:)_{n^2,1} \in \CC^{n^2} \ . $\\[0.5mm]
\hspace*{5mm} \ With $P(t)  = (X^T(t) \otimes I_n + I_n \otimes X(t)) \in \CC_{n^2,n^2}$ assumed nonsingular for all $t$ and\\[0.5mm]
\hspace*{5mm} \ $ q(t) = \dot A(t)(:) + \et\ A(t)(:) - \et\ (X^T(t) \otimes I_n) \cdot X(t)(:)  \in \CC^{n^2}$ we finally have separated \\[0.5mm]
\circledthree$\!\!_{\small nonsing}$ \hspace*{40mm}  $\dot X(t) = P(t)\backslash q(t)$.\\[0.5mm]
\hspace*{5mm} \ If $P(t)$ is singular for some $t$, we cannot invert $P(t)$  and need to use $P$'s pseudoinverse instead. For more\linebreak
\hspace*{5mm} \  details on this process, see \cite[Section 3, (III), (IV), (VII) and (VII start-up)]{FUsurveyZNN}.\\
\circledfour \ Choose a convergent look-ahead finite difference formula of type {\tt j\_s} that involves $X(t_{k+1}), X(t_{k-1}), X(t_{k-1}), ...$ \linebreak
\hspace*{5mm} \ and $\dot X(t_k)$.\\

\circledfive \ Solve the difference formula in \circledfour algebraically for $\dot X(t_k)$.\\
\circledsix \ Equate both expressions of $\dot X(t_k)$ in steps 
\circledthree$\!\!_{\small nonsing}$ and \circledfive and solve the resulting solution-derivative free\linebreak
\hspace*{5mm} \  equation for $X(t_{k+1})$ numerically by solving the linear system in \circledthree$\!\!_{\small nonsing}$ and adding the recursion in \circledfive$\!\!$.\\
\circledseven \ Iterate step \circledsix for $t_{k+2}, t_{k+3}, ... $ until you reach the final time $tf$.\\[2mm]
Discretized ZNN benefits from two important innovations. The first is the exponential error decay stipulation and the setting of $\tau$  in \circledtwo$\!\!$. The second are the combined workings of \circledthree through \circledfive where ZNN places a barrier,  a membrane between the derivative based formulation and their conditioning problems and the linear equations solving plus recursion step \circledsix$\!\!$. Most simply said: conditioning problems at the error equation  or the error equation ODE level are not passed into or onto the actual solving part of ZNN and do not affect the linear system conditioning in \circledsix$\!\!$. The two -- quite different  -- conditioning limitations of ODE solvers and of linear systems' solvability  are not shared and thereby ZNN's output and efficiency  is immune to original problem instabilities and stiffnesses et cetera in the initial problem set-up  \circledzero$\!\!$, \circledone$\!\!$, and of the ODE set-up in  \circledtwo$\!\!$, \circledthree$\!\!$.\\[2mm]
Our first desire was to learn how to assign separate values to the decay constants $\eta$ used in the start-up phase and in the regular ZNN iteration phase. After some experiments and trial and error discoveries, we have settled on two different exponential decay setting of $\eta$, namely  at 160 for the start-up phase that we  run for 12 steps with Euler while our chosen finite difference formula of type {\tt j\_s = 4\_5} requires  only $j + s = 9$ actual solution data points. We have adjusted $\eta$ in the iteration phase of ZNN to $\eta = 1.45$ which  differs slightly from the earlier global setting of $\eta = 1.35$ for Figure 1 left. Our time interval for Figure 3 is 1 hour or 3600 seconds. The error drops below $10^{-11}$ and stays below after 23 seconds and levels out around $10^{-13}$ at the end of the run.
Our simulation run took less than 14 seconds and it computed  18,000 points of the solution $X(t_k)$ using around $8 \cdot 10^{-4}$ seconds per time step $t_k \in [10,3610]$, or in less than 4 \% of the constant sampling gap of $tau = 0.02$ that simulates a discrete 50 Hz sensor input.\\
The error function graphs of Figure 1 left and Figure 3 look similar  but the AZNN method reaches a ten fold lower error plateau  when compared to basic ZNN. In Figure 3 AZNN's $\eta$ values were altered between phases and the number of initial Euler steps was increased from the necessary start-up data length $j + s = 4 + 5 = 9$ to 12 for this time-varying matrix flow problem.\\[2mm] 
We start Euler from the entry-wise square root matrix {\tt 
\verb*/A(to).^{0.5}/}  in Matlab notation for ease with on-chip implementation.\\[-8mm]
\begin{center}
\includegraphics[width=112mm]{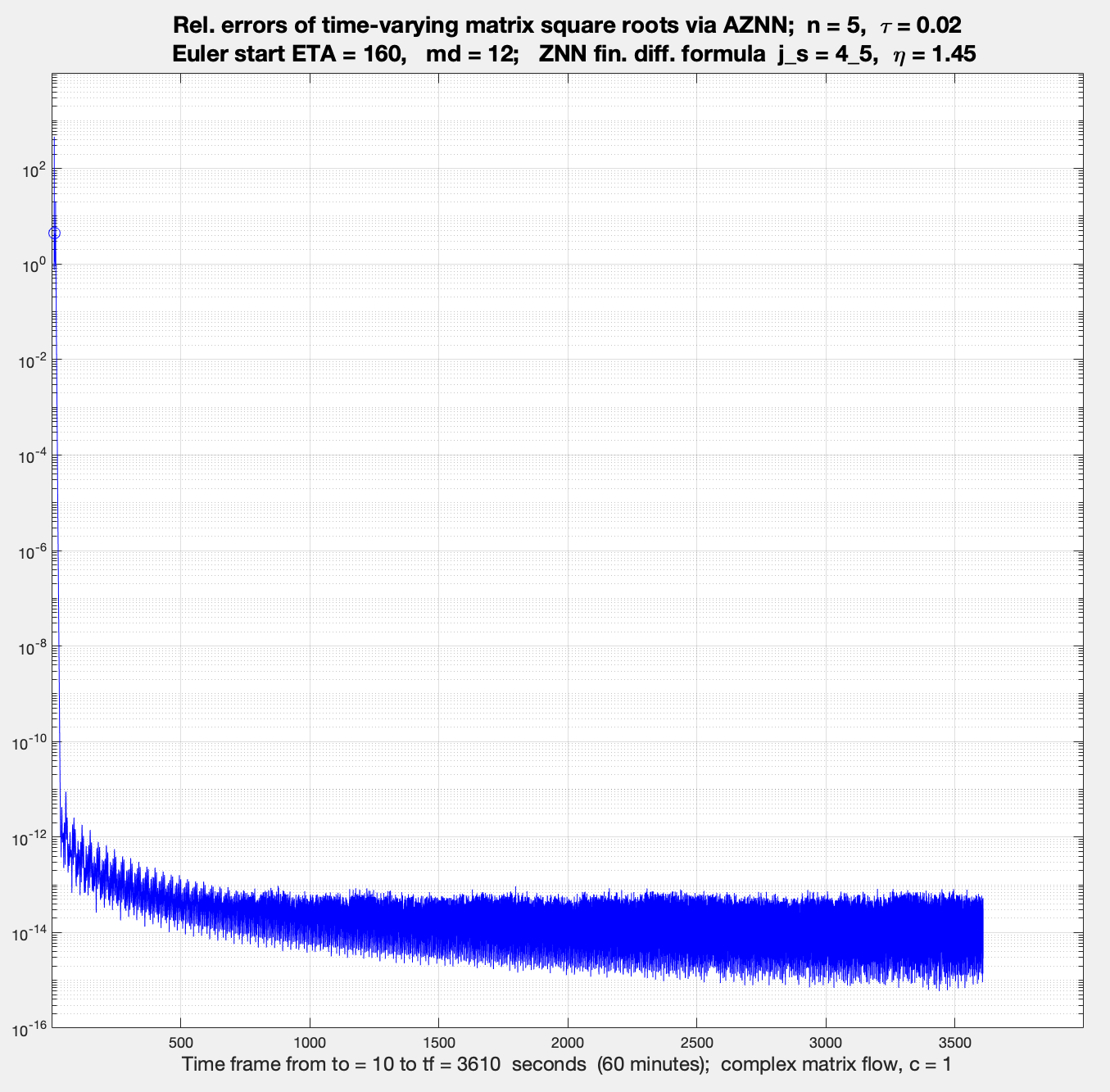}\\[-1mm]
Figure  3 
\end{center}

All of our finite difference formulas are convergent and as convergent difference formulas they obey certain rules. For example with the 5-IFD formula $$ \dot z_k = \dfrac{8z_{k+1}+z_k -6z_{k-1} - 5 z_{k-2} + 2 z_{k-3}}{18 \tau} \in \CC^{n+1} $$
of type {\tt 2\_3}, its global truncation order $O(\tau^3)$ and its scaled numerator above expressed as\\[-4mm] 
$$ z(t_{k+1}) + \al_k z(t_k) + \al_{k-1} z(t_{k-1}) + \ \cdots \ + \al_{k-\ell} z(t_{k-\ell}), $$

\vspace*{-1mm}
has the characteristic polynomial $ p(x) = x^{k+1} + \al_k x^k + \al_{k-1} x^{k-1} + \ \cdots \ + \al_{k-\ell} x^{k-\ell} $. Convergence now requires that $p(1) = 1 + \al_k + \al_{k-1} + ... + \al_{k-\ell} = 0$, see \cite[Section 17]{EMU96}. Plugging $z_{..}$ into the 5-IFD difference formula, we realize that    \\[-2mm]
$$\hspace*{34mm}z(t_{k+1} ) +\dfrac{1}{8} z(t_k) - \dfrac{3}{4}  z(t_{k-1}) - \dfrac{5}{8} z(t_{k-2})  + \dfrac{1}{4} z(t_{k-3}) \approx o_{n+1} \hspace*{30.8mm} (*)$$
is approximately the zero vector due to unavoidable truncation and rounding errors.
The recursion part of the ZNN computations is split by an equal sign in the final computational equation in step \circledsix. Thus $$|z(t_{k+1} )| \approx  \left|-\dfrac{1}{8} z(t_k) + \dfrac{3}{4}  z(t_{k-1}) + \dfrac{5}{8} z(t_{k-2})  - \dfrac{1}{4} z(t_{k-3})\right| .$$
And the two nearly equal recursion parts above  become ever closer to each other in magnitude due to the convergence to zero condition $(*)$ of all convergent finite difference schemes. Thus the first linear equations solution term  $ (P(t_k) \backslash q(t_k))$ in \circledsix eventually has only a  comparatively small magnitude.  It adjusts the predicted value of $z(t_{k+1})$ only slightly  according to the current system data inputs while the remaining finite  difference formula term $$ -\dfrac{1}{8} z_k + \dfrac{3}{4} z_{k-1} + \dfrac{5}{8} z_{k-2} - \dfrac{1}{4} z_{k-3}$$  has nearly the same magnitude as  the computed solution vector $z_{k+1}$ for large indexed time steps  $t_{k+1}$.\\[2mm]
This heuristic analysis explains the small but persistent error magnitude  wiggles that set in after a few dozen seconds in Figures 1 and 3 above when the error level has dropped to well below $10^{-10}$. Towards the end of our  hour long simulation, the linear equations term that is computed in the AZNN iterative matrix square root example  for $z_{k+1}$ has only about 1/10,000 of the magnitude of the finite difference formula  term while our trial matrix flow norms have grown to around $10^6$.  The linear equations  'correction term' in $z_{k+1}$  then affects the average number of 13 accurate digits only slightly. These variations by around 1 digit up and down are depicted in Figures 1 and 3  as 'wiggles' and they assure us of between 12 and 14 accurate digits   in all sufficiently late computed solutions.\\[1mm]
Following Figure 4 below for a matrix symmetrizer problem we include a data table  that describes the same theoretical magnitude disparity.
Our Matlab codes for AZNN computations of time-varying matrix square roots are\\ tvMatrSquareRootwEulerStartva.m and AZNNtvMatrSquRootEuler5.m.\
 Both are available  at \cite{FUMatlabAZNN} together with the auxiliary  codes  Polyksrestcoeff3.m, formAnnm.m, and formAnnmc.m.\\[1mm]

{\large \bf (B) \ \ \   Time-varying symmetrizers $S(t)$ for time-varying matrix flows $A(t)$ via AZNN}\\[2mm]
Matrix flow symmetrizers are much easier to handle in AZNN as shown  in Figures 4 and 5 below for example. Here the separate exponential decay parameters of AZNN can apparently be set rather high, or even astronomically high when compared with the rather moderate decay constants that are suitable for  convergence  in the same  matrix flow example $A(t)$  in Section 2 (A) for finding time-varying matrix square roots in Figure 3. The data or linear solver term induced wiggles of $z(t_{k+1})$ are much smaller and less noticeable here. They occur in a relatively small time interval from about 2 minutes to about 45 minutes after start-up. In this time interval, the magnitudes of the linear solver term in the solution remains around 4 digits below that of the finite difference formula's  term. After around 6 hours the magnitude disparity of the two terms increases to 5 digits making the contribution of the linear equations term solve almost irrelevant. But the exponential truncation errors decrease is much slower even  with ten or hundred times larger decay constants than those that we found for the time-varying matrix square roots problem.\\[1mm]
 Individual matrix problems for the same data may behave quite differently under AZNN.\enlargethispage{20mm}

% using  ... 
%   Box/matlabin/AZNNadjusted/AZNNtvMatrSymmwEulerStartm5
%  in Section 2 (B)
\vspace*{-2mm}
\begin{center}
\includegraphics[width=118mm]{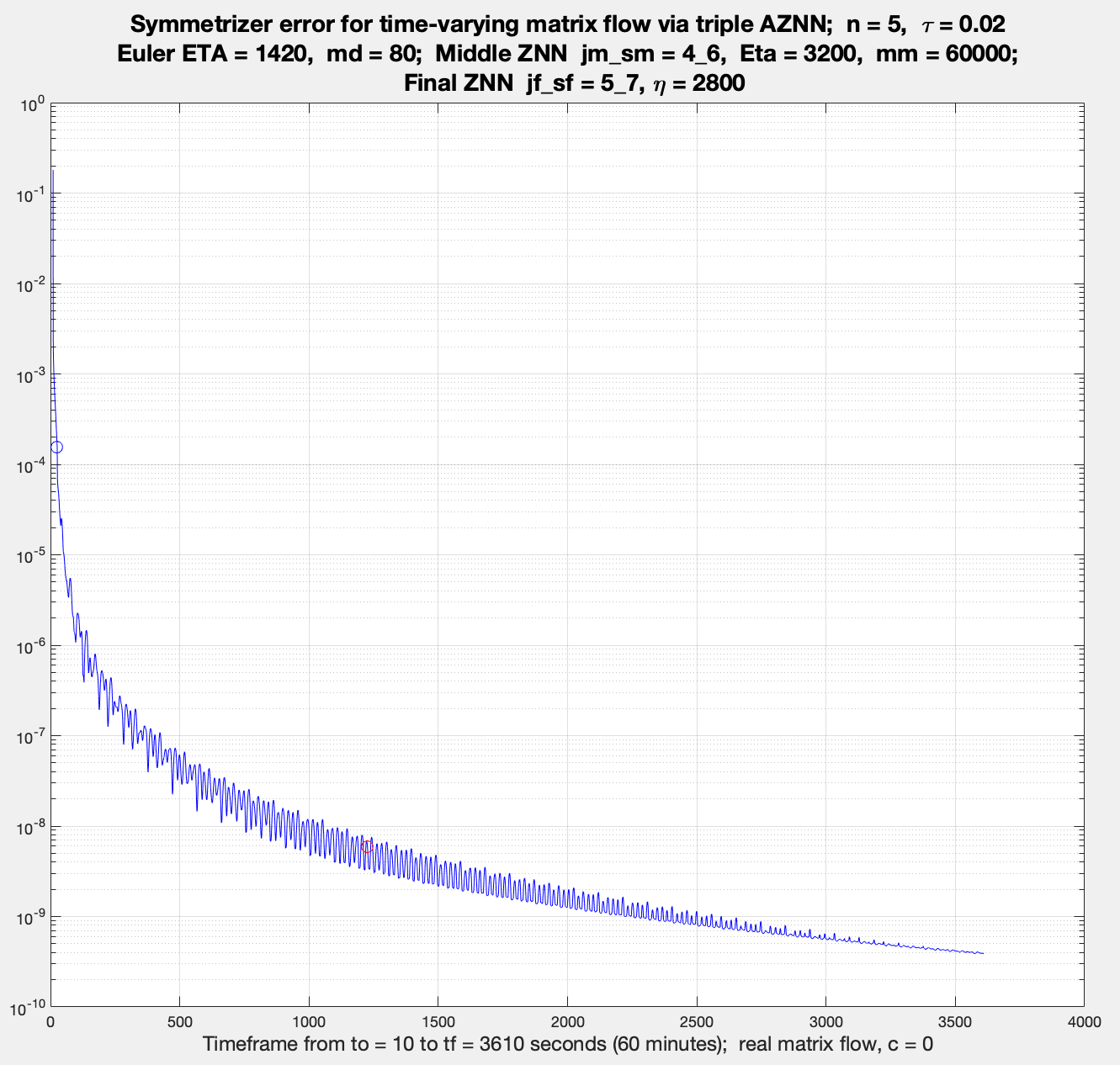} \\[-0mm]
Figure  4 \end{center}

\vspace*{-3mm}
Here is a table of convergence data and magnitude comparisons for the linear equations solution part and the partial difference recurrence part that drive the iterations  for the symmetrizer example  in Figures 4 and 5.\\[-12mm]

\setlength{\unitlength}{1mm}

\newpage

\vspace*{-7mm}
\begin{center}
\includegraphics[width=142mm]{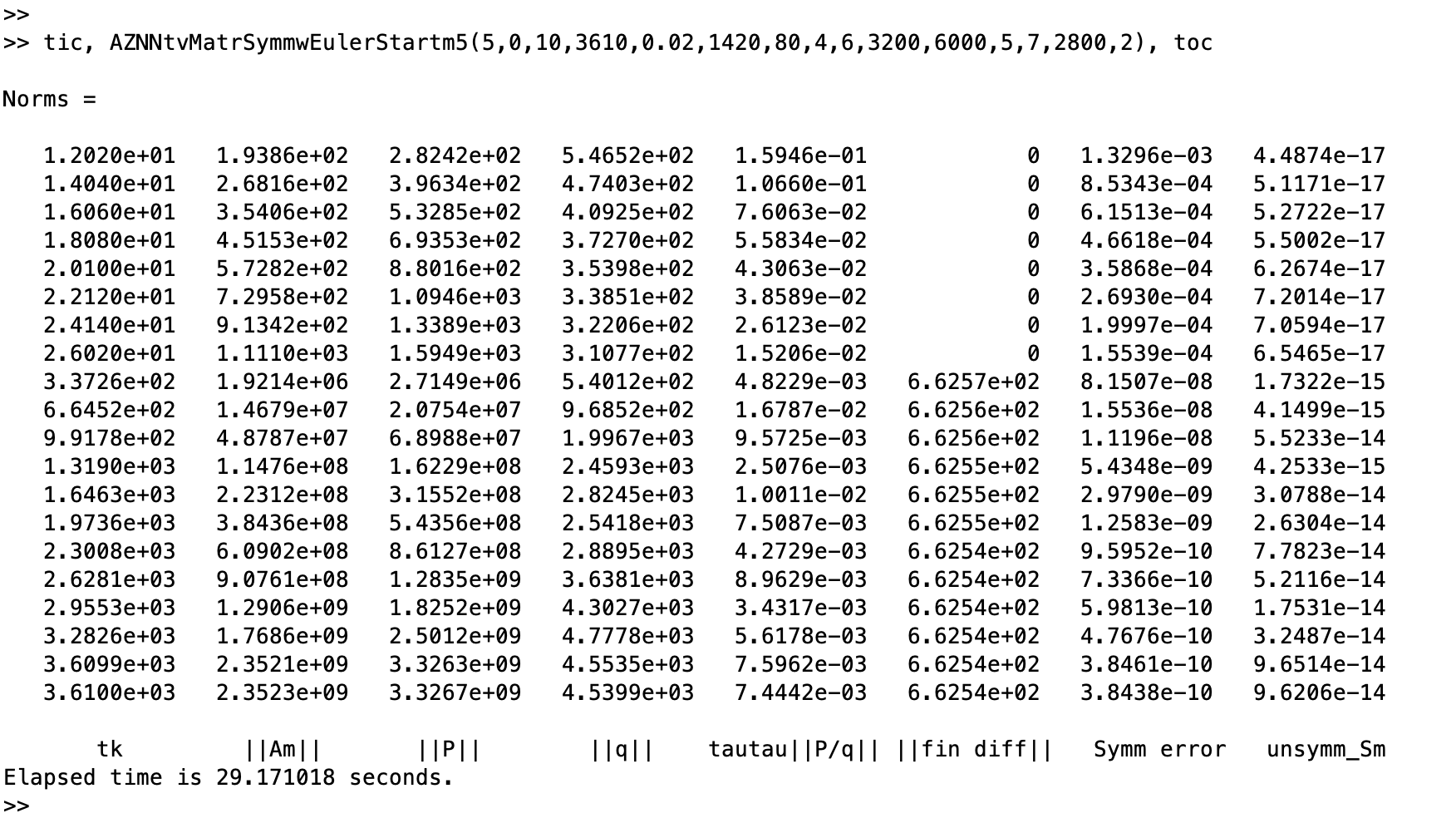} \\[-2mm]
Table  1 \end{center}

\begin{picture}(80,100)(0,-111)

\put(77.2,0){\framebox(35,62)[t]{ }}

\end{picture}

\vspace*{-103mm}
The columns in Table 1 list computational data for several discrete times $t_k$. The second column contains the matrix norms of the time varying input matrices $A(t_k)$, followed by the norm of the linear equations' system matrices $P$  and of their right hand sides $q(t_k)$ when computing the next solution at $t_{k+1}$. The two framed columns 5 and 6  contain the magnitudes of the two terms that predict the symmetrizer $S(t_{k+1})$ of $A(t_{k+1})$. The final two columns measure the actual symmetrization error $|| A(t_{k+1}) \cdot S(t_{k+1}) -  S^T(t_{k+1}) \cdot A^T(t_{k+1})||$ for $A(t_{k+1})$  and the un-symmetry of $S(t_{k+1})$. Clearly the two additive terms in the ZNN formulation of  $S(t_{k+1})$ have very different magnitudes and the finite difference term is approximately 10,000 times larger than the linear equations term. This data corroborates  our earlier theoretical analysis very well.

\vspace*{-1mm}
\begin{center}
\includegraphics[width=99mm]{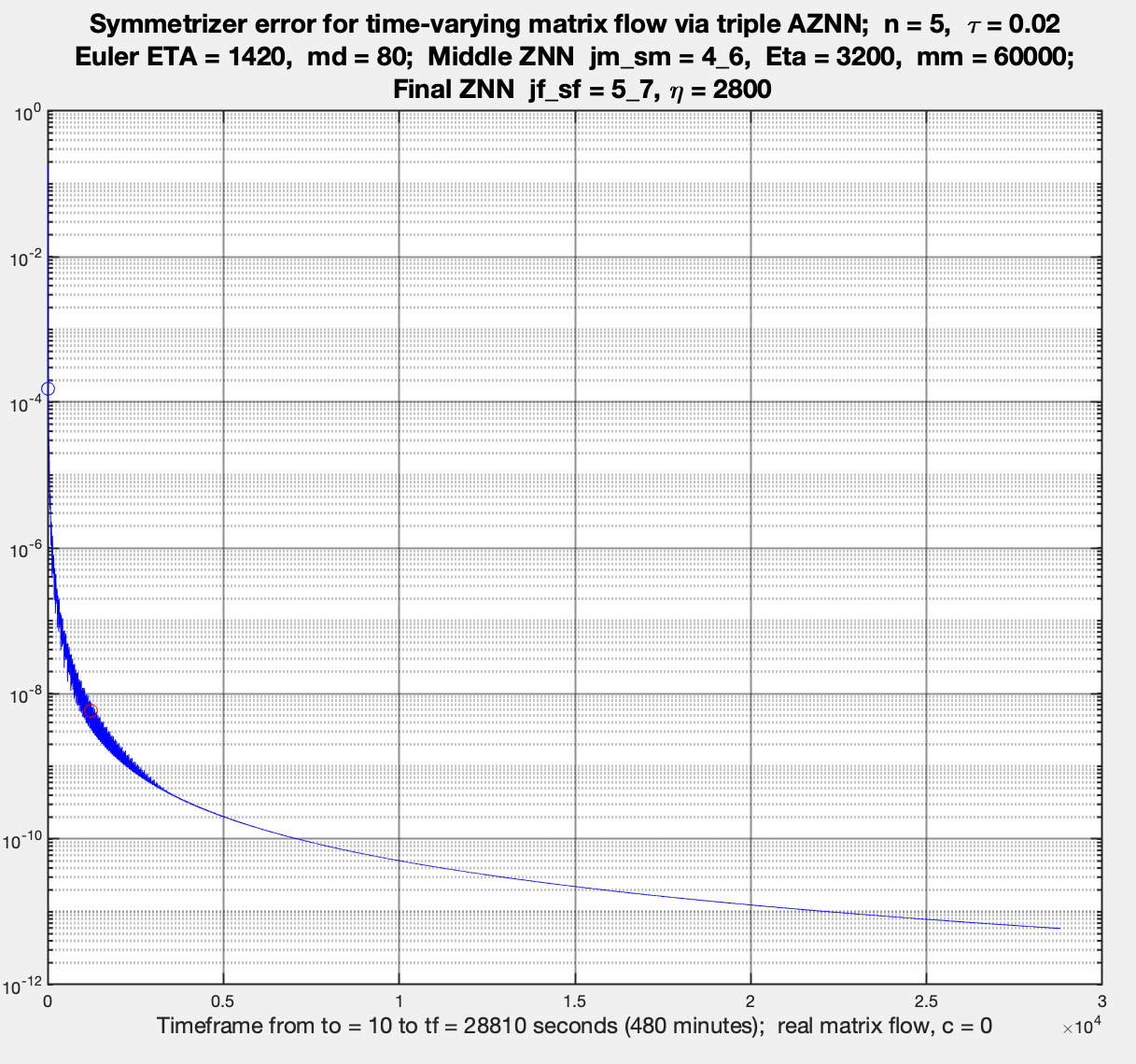} \\[-0mm]
Figure  5 \end{center}

\vspace*{-2mm}
Our next example deals with finding a complex symmetric matrix factorization $S(t) \cdot V(t)$ of a complex matrix flow $A(t)$ via AZNN in a 'wiggle-free' way if possible.\\[1mm]
Here we use one fixed complex matrix flow $A(t)$ and determine the $\eta$ intervals that ensure wiggle-free truncation error curves over long time intervals. In our tests  with AZNN below we keep the Euler initial part value of $\eta$  and the length of the start-up iterations fixed  while we vary the value of $\eta$ in the subsequent basic ZNN iterations.\\[1mm]
% use  ...   ../Box/...matlabin/AZNNadjusted/AZNNtvMatrSymmwEulerStartm4   below
\hspace*{-1.5mm}\includegraphics[height=71mm]{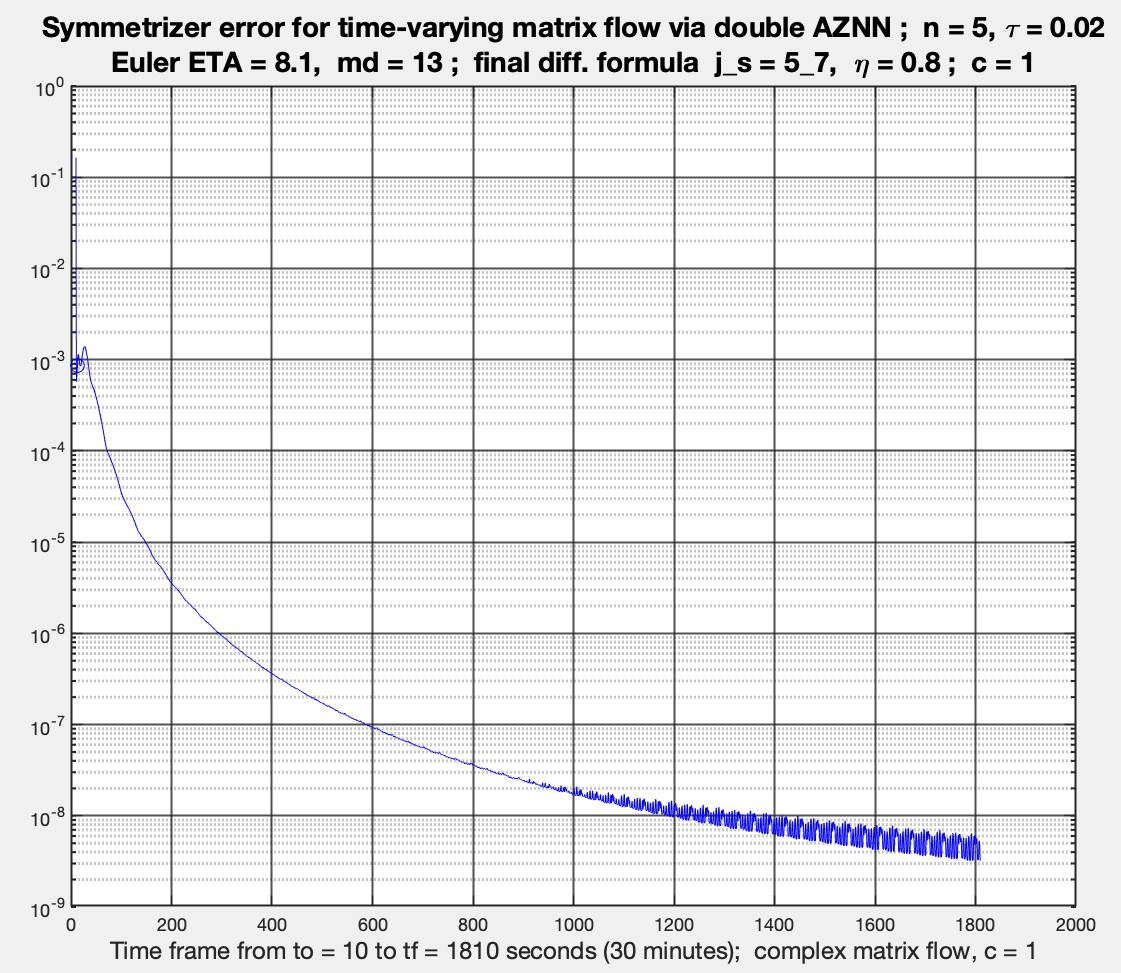} \vspace*{-2mm}\includegraphics[height=71mm]{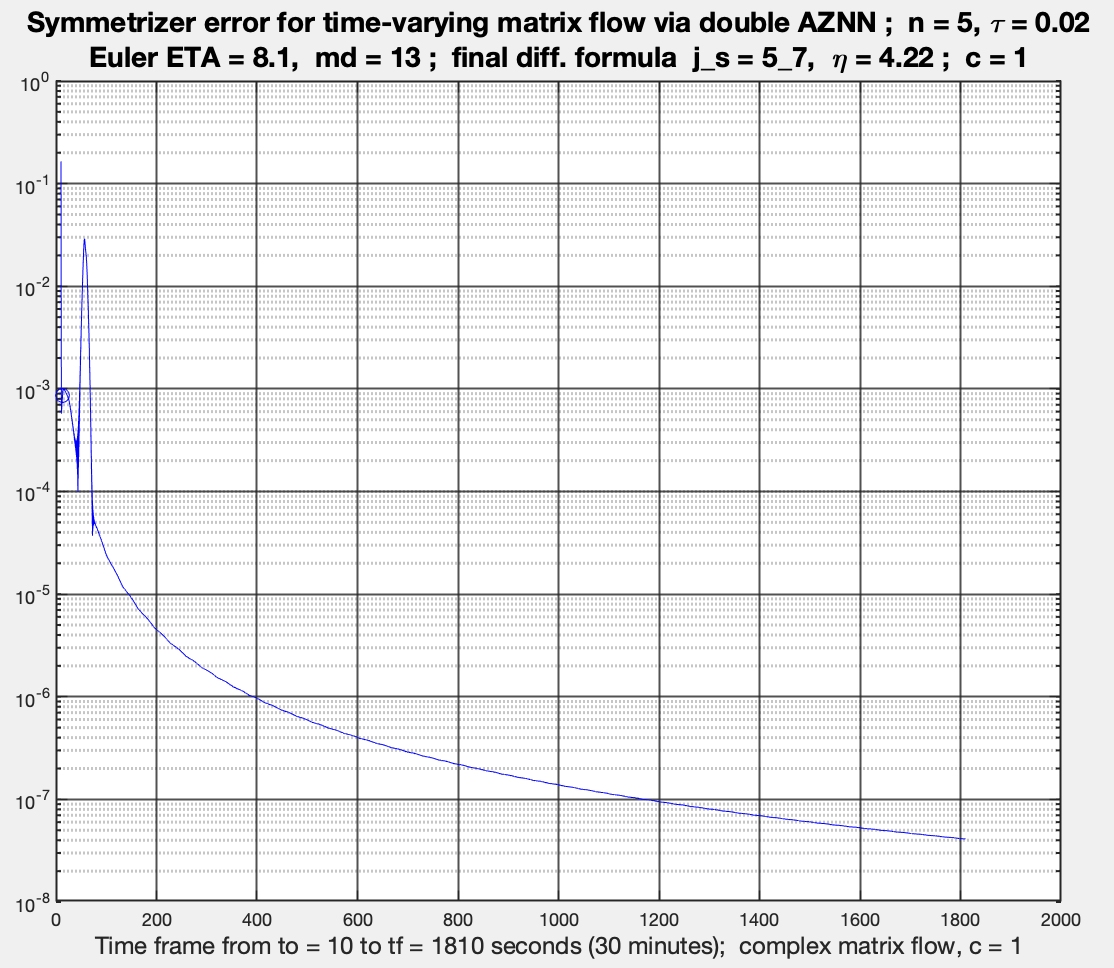}
\\[1mm]
\hspace*{73mm} Figure  6\\[1mm]
Figure 6 depicts the extreme lower and upper values for $\eta$ that induce wiggly truncation error curves for this example flow. Next we  show upper and lower bounds for non-wiggly truncation error curves, except possibly right near the Euler to ZNN transit point in time,  in between which any $\eta$ can be chosen for smooth output.\\[1mm]
% use  ...   ../Box/...matlabin/AZNNadjusted/AZNNtvMatrSymmwEulerStartm4   below
\hspace*{-1.5mm}\includegraphics[height=71mm]{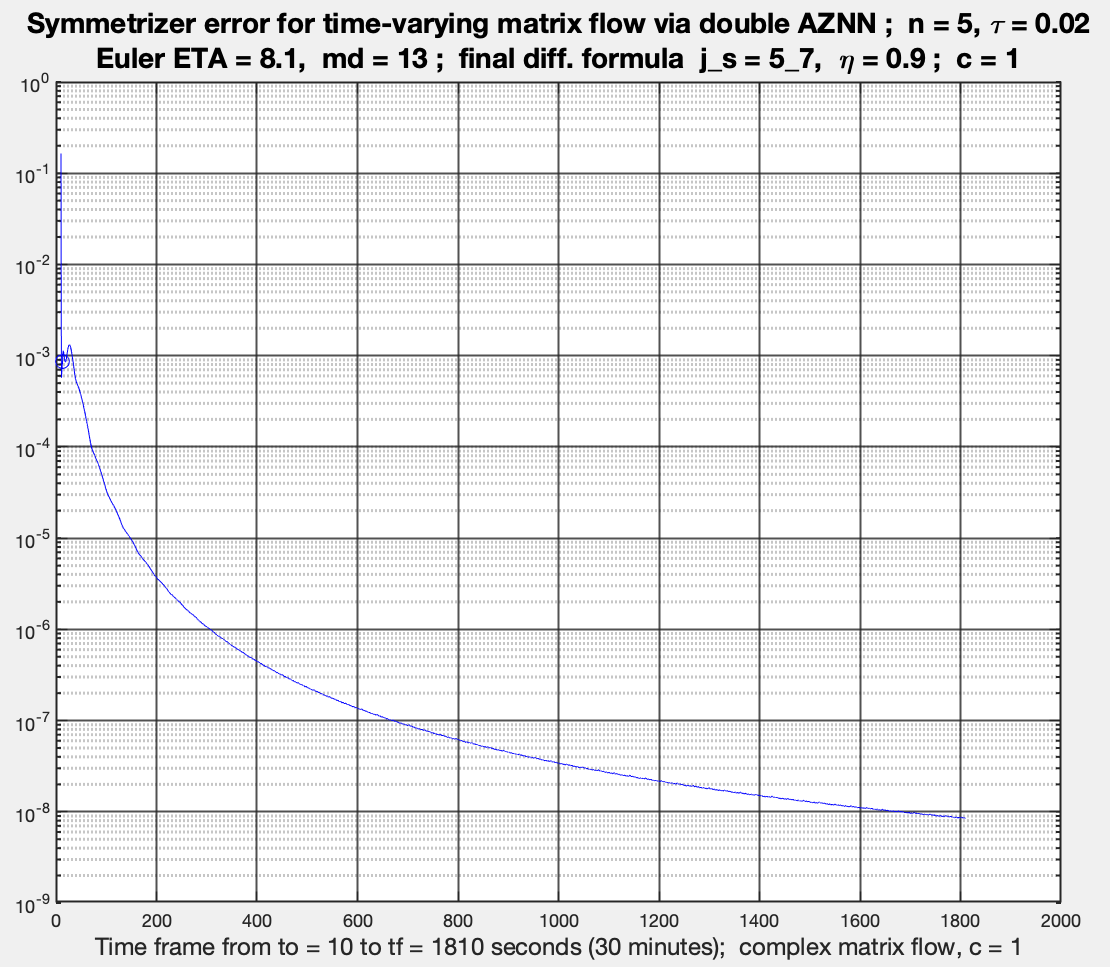} \vspace*{-2mm}\includegraphics[height=71mm]{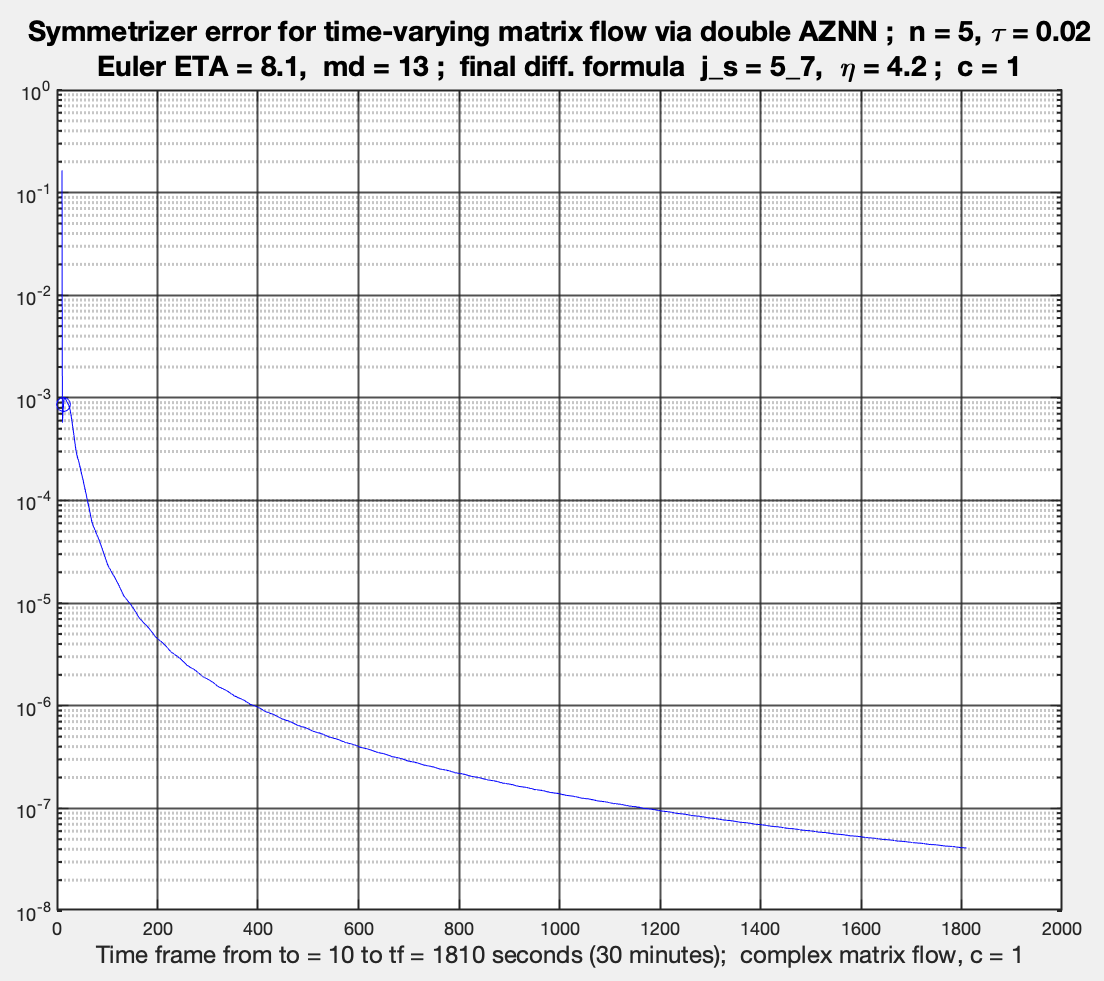}
\\[1mm]
\hspace*{73mm} Figure  7\\[1mm]
In Figure 7 the lowest relative symmetrizer  errors are recorded in the left graph as $8.4872 \cdot 10^{-9}$ at the end of the 30 minute run when $\eta = 0.9$ and as $4.1151 \cdot 10^{-8}$ in the right graph for $\eta = 4.2$. Here the slower exponential decay constant gives us a five times  better convergence result. This seems to turn our intuitively expected result on its head. And we know of no explanation for this mystery.\\
For  5 by 5 real matrix flows the AZNN computations of the flow's  symmetrizers take around 28 seconds per 1 hour of simulation or  224 seconds for the 8 hour plot in Figure 8.  It  reached  a minimum relative
symmetrization error of   $5.9095\cdot 10^{-8}$ a the end of the run.\\[-5mm]
\begin{center}
\includegraphics[width=108mm]{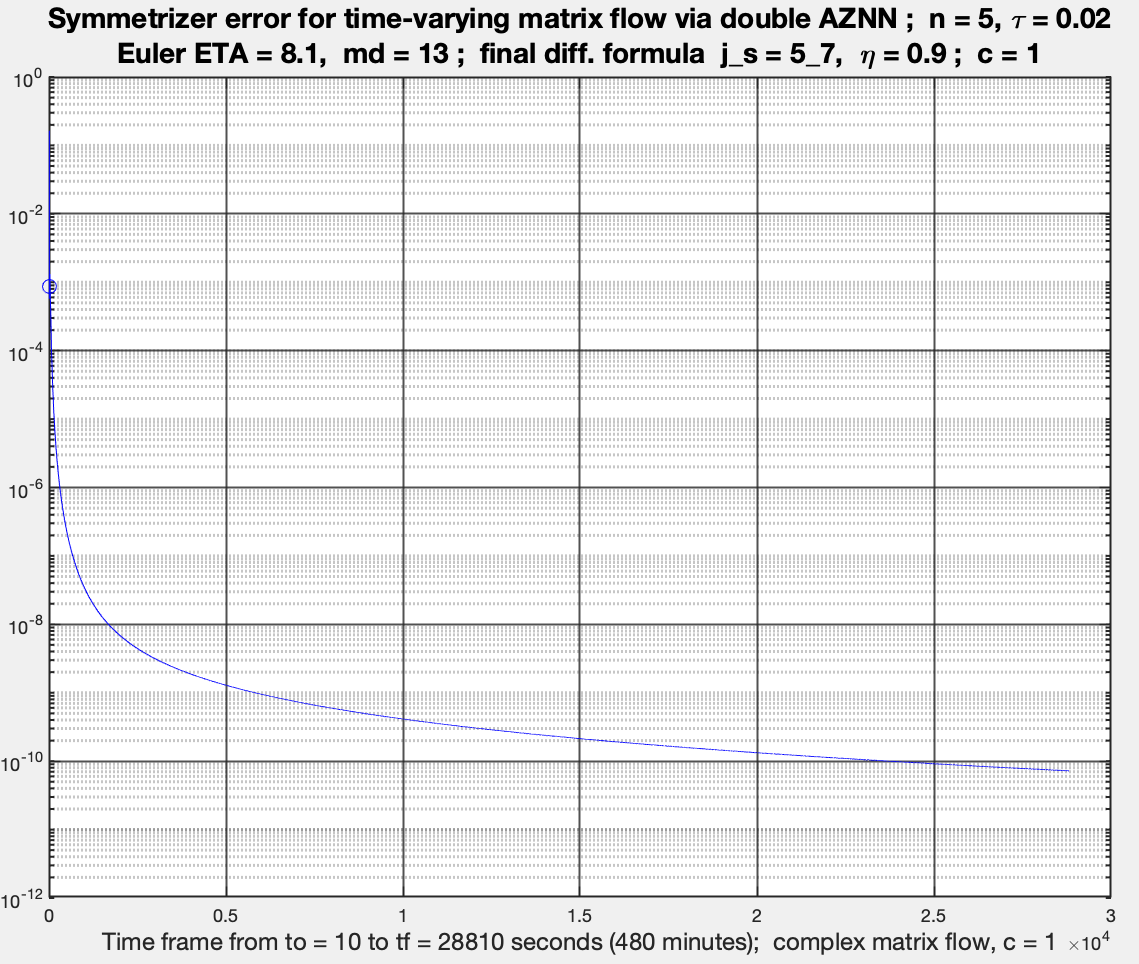} \\[-1mm]
Figure  8 \end{center}

\vspace*{-2mm}
Our 5 by 5 complex matrix flow experiment for time-varying matrix symmetrizers  takes around 48 seconds per 1 hour of simulation or  385 seconds for the 8 hour plot in Figure 8  that  reached  a minimum relative
symmetrization error of   $7.2301\cdot 10^{-11}$ a the end of its run.\\[2mm]
Exploring some phenomena of AZNN in this paper has lead  to several new open problems with neural network based time-varying matrix algorithms:  parameter optimizations in AZZN appear to be problem dependent; for any one problem large magnitude variations are seemingly possible for decay constant  intervals with comparative results, even thousand-fold $\eta$  variations; why are there wiggles in some relative AZNN error graphs and so forth.\\[2mm]
Matlab codes for AZNN computations of time-vaying matrix symmetrizers are AZNNtvMatrSymmwEulerStartm5.m for triple phase parameter settings and AZNNtvMatrSymmwEulerStartm4.m for double parameter settings for the start-up Euler phase and the ZNN iterations phase only. They are available together with trial matrix flow  generating codes  at \cite{FUMatlabAZNN}.\\[1mm]
 A new  application of AZNN methods uses time-varying neural networks  directly to solve some previously unsolvable static matrix symmetrizer problems in Section 3.\\[-7mm]

\section {Computing static matrix symmetrizers via time-varying AZNN}

We start with general observations of matrix factorizing methods such as LR, QR, or SV.\\[1mm]   
 For a given fixed entry matrix $A_{n,m}$, real or complex,  almost every  matrix factorization  $A = X \cdot Y$ with specified qualities in $X$ and $Y$  is achieved by operating on the given matrix $A$ from the left with a matrix $U$ of desired type to obtain  the right factor $Y$ as $U \cdot A = Y$ with $Y$ as desired and then inverting $U$ to obtain $A = X \cdot Y$ with $X = U^{-1}$ if the inversion of $U$ preserves its special type. To factor $A = L \cdot R$ for lower and upper triangular matrices $L$ and $R$ for example, we obtain the upper triangular $R$ by repeatedly multiplying $A$ and its updates by lower triangular Gaussian elimination matrices from the left. Their product is again lower triangular and it can be easily inverted while keeping the lower triangular form to give us  the LR factorization of $A$. Likewise for  orthogonal Householder matrices that reduce $A_{n,m}$ from the left to upper triangular form in the QR matrix factorization.\\
 The key issue for successfully expressing any given real or complex square matrix $A$ as a product of two symmetric matrices $S$ and $V$ as $A =  S \cdot V$ is to find a nonsingular symmetric matrix $U = U^T$  with $U \cdot A = V = V^T = A^T \cdot U$. Then  a symmetric matrix factorization of $A$ is  $A = S \cdot V$ with $S = U^{-1}$. In 1902  Frobenius \cite{F1910} proved theoretically that such a factorization exists  for each square matrix $A$. But finding {\em nonsingular} left side  symmetrizers   numerically  for all static matrices $A_{n,n}$  proved quite elusive, see  \cite[p. 609-612]{DU16}, especially for defective or derogatory matrices $A$ with repeated eigenvalues and certain Jordan structured ones inside a recent extensive comparison of  eight iterative and eigen or svd based symmetrizer algorithms \cite{DU16}.\\[2mm]
 An application of  time-varying AZNN to compute matrix symmetrizers for a given fixed entry matrix $A$ is relatively simple. We use AZNN for the associated matrix flow ${\tt A(t) = t*A + (1-t)^a*BB}$ here with an exponent $a > 0$ that we choose heuristically through experimentation with $A(t)$ and a random small entry matrix ${\tt BB}$.  We start the AZNN process with Euler and a chosen constant sampling gap $\tau   = t_{k+1} - t_k$ from  $t_o < 1$. Typically $t_o = 0.985$ when using small $\eta$ and $\tau$ parameters and $t_o = 0.9985$ for large $\eta$ and $\tau$ parameters. We stop the AZNN iterations at $t_k = 1$ because then we have computed a symmetrizer $S$ of $A(1)$ via AZNN. Clearly this symmetrizer also symmetrizes $A$ since $A(1) = A$ by construction. Here it is important to synchronize $t_o$ and $\tau $ so that in the iteration process ${\tt A(t_k) = t_k*A + (1-t_k)^a*BB}$,  the factor $1-t_k$ reaches  zero precisely for some index $k$.\\[1mm]
 Now we test AZNN for the  fixed entry matrix examples in \cite{DU16} for which none of its eight symmetrizer algorithms could find a nonsingular symmetrizer. The first one is the Kahan matrix of dimension 35, see  \cite[p.609, first table]{DU16} and the {\em gallery} in Matlab.\\[1.5mm]
\hspace*{-1.5mm}\includegraphics[height=68mm]{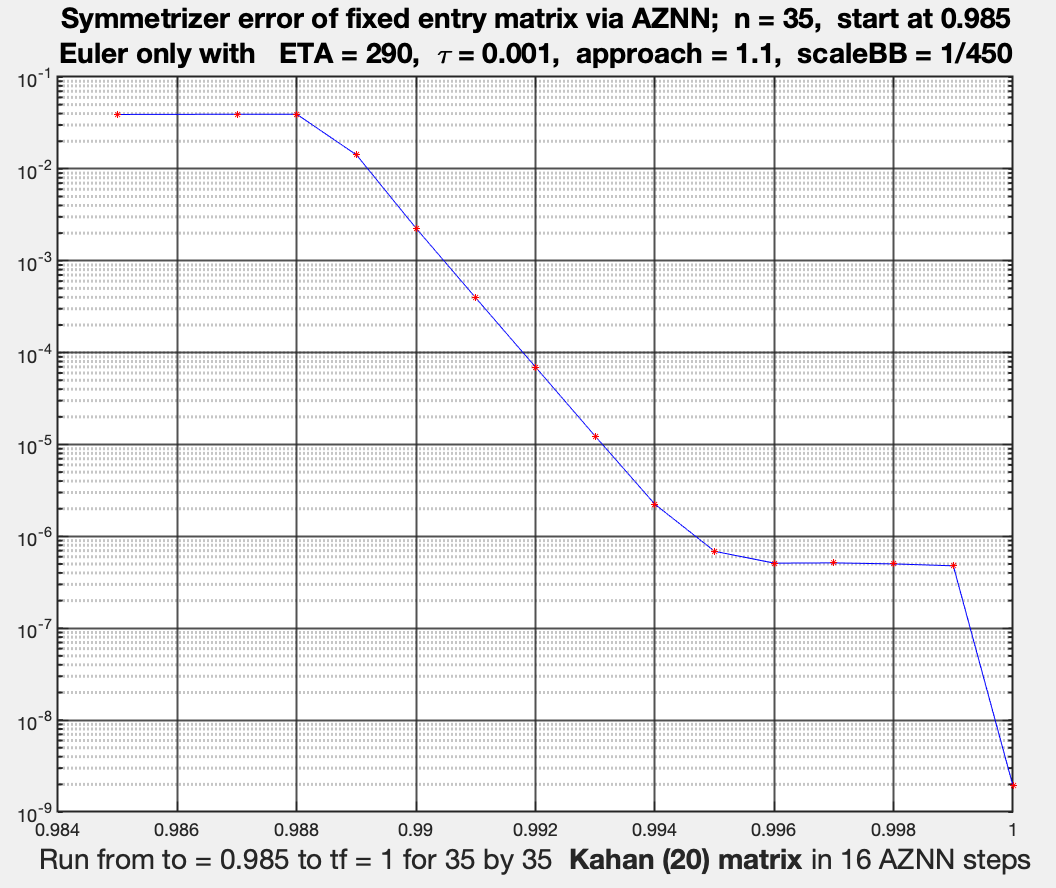} \vspace*{-2mm}\includegraphics[height=68mm]{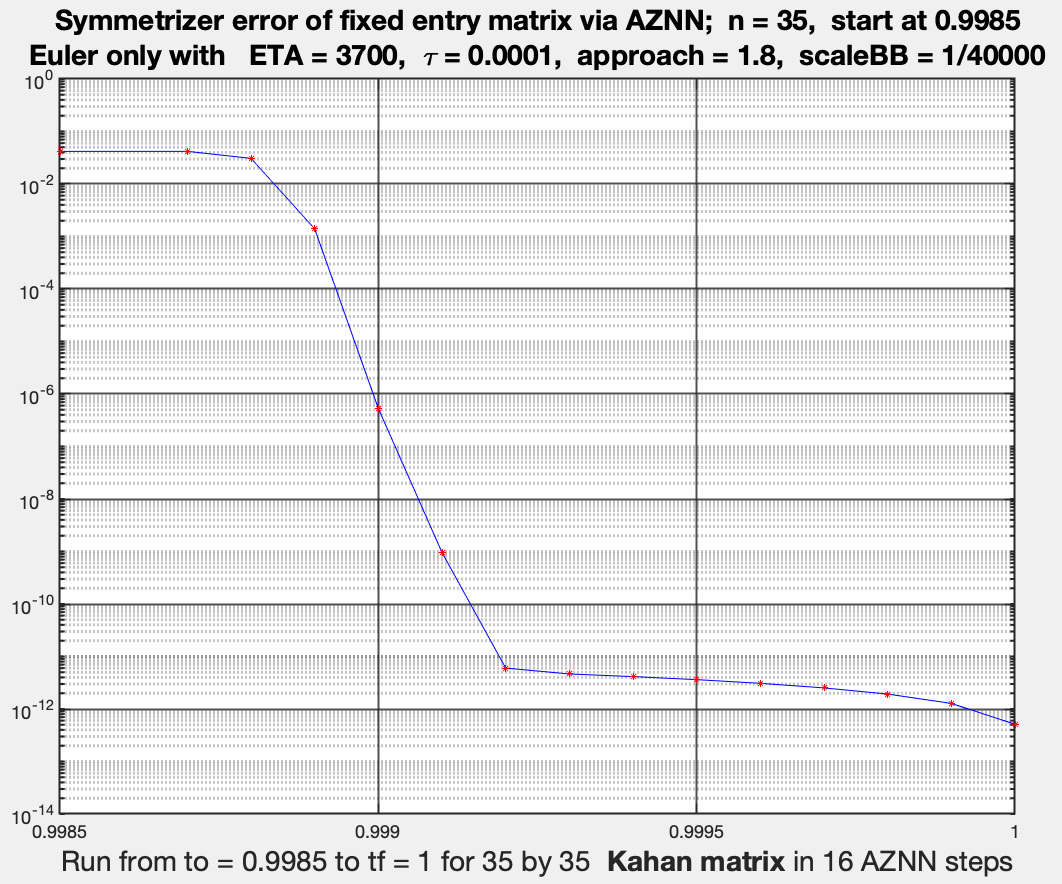}
\\[1.5mm]
\hspace*{73mm} Figure  9\\[1.5mm]
The left plot of Figure 9 shows the relative error  $||S(1)\cdot A - A^T\cdot S(1)||/||A||$ for the computed Kahan (20)  matrix symmetrizer $S = S(1)$ at the end point  $t_k = 1$ of AZNN for a relatively small exponential decay parameter ETA and a relatively large first ${\tt BB}$ perturbation for $A(0,985)$. The value of the {\em approach exponent} $a$ in $(1-t_k)^a$ in the definition of $A(t_k)$ here is set to 1.1 for all iterates, indicating a convex, slightly quadratic function approach to $t_k = 1$ from below a linear approach with $a=1$. Sixteen AZNN steps take us into the $10^{-9}$ error territory here. Our computed symmetrizers $S(t_k)$ for 35 by 35 Kahan all have condition numbers around $10^8$ or $ 10^{10}$ and full rank, while the ones computed in \cite[p. 609]{DU16}   in six different static matrix ways have condition numbers $10^{9}$ or $10^{10}$ with full numerical rank according to Matlab's criteria and four eigen method based solutions have  symmetrizer condition numbers  of $10^{15}$  and deficient numerical rank 32 only. The 2-norm of the 35 dimensional Kahan matrix is 4.8.\\
The right side plot with increased parameters for Eta and the approach constant  $a$, both  by  factors of around 10, and with a 100 fold decreased  perturbation {\tt BB} size and a 10 fold smaller sampling gap $\tau$, i.e., with a countervalent change of parameter pairs,  brings the relative error for AZNN down into $10^{-13}$ territory and computes a full rank symmetrizer for the Kahan matrix.\\[1mm]
Next in \cite[p.409, second table]{DU16} comes the Frank matrix of size 35 by 35 and 2-norm 340. The plots in Figure 10 show quite similar error behavior as the Kahan matrix did with or without countervalently set parameters, except for the much smaller relative symmetrizing errors that fall well below the machine constant in the right  graph of Figure 10. \\[1.5mm]
\hspace*{-1.5mm}\includegraphics[height=68mm]{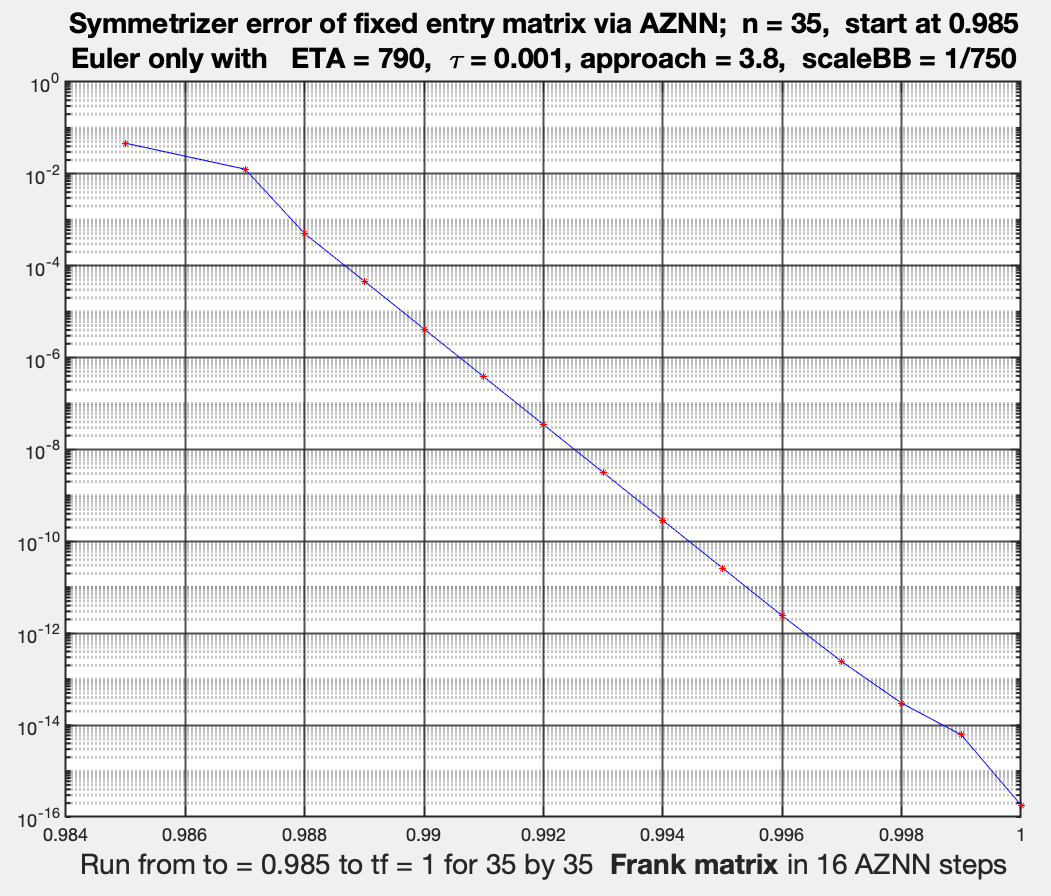} \vspace*{-2mm}\includegraphics[height=68mm]{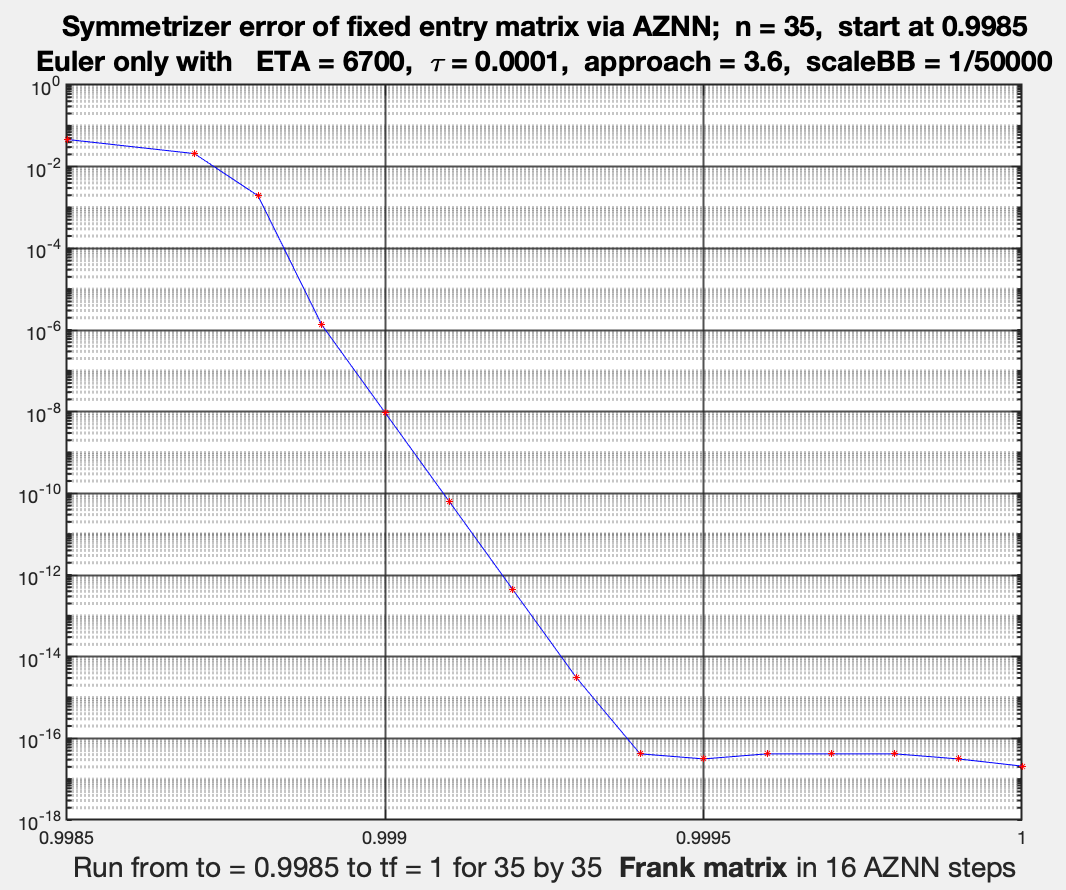}
\\[1.5mm]
\hspace*{73mm} Figure  10\\[1.5mm]
A derogatory 23 by 23 upper triangular matrix $C_{ut}$ with 2-norm 281 was studied in \cite{DU16}. Here are our AZNN results for this static matrix symetrizer problem  with our standard large and small parameter value pairs.\\[-5mm]
\begin{center}
 \includegraphics[height=80mm]{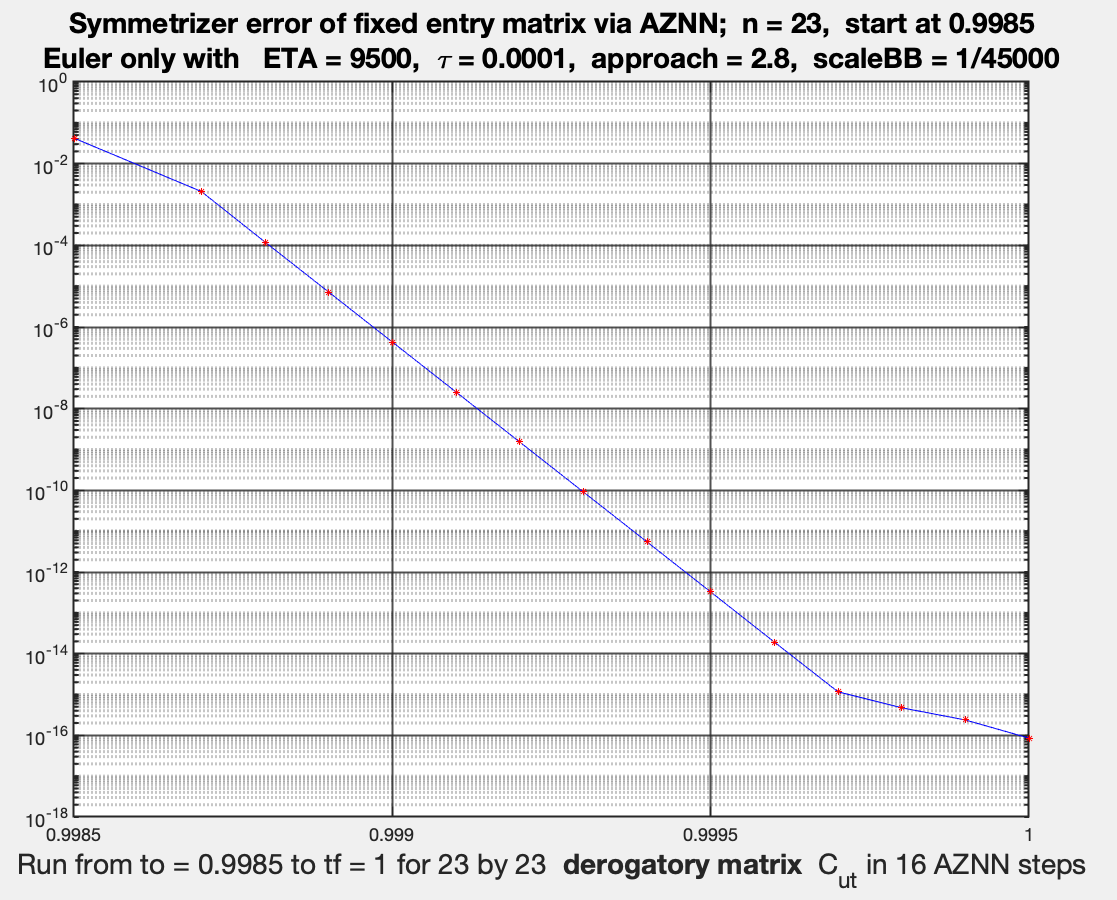}\\[-1mm] 
 Figure  11
\end{center}

\vspace*{-1.5mm}
In Figure 11 we have again used relatively large and small countervalent parameter pairs together to compute a nonsingular symmetrizer of $C_{ut}$ with symmetrizing error below the machine constant and condition number 362. In \cite[p. 410, second table]{DU16} all six iterative and eigen based methods compute $C_{ut}$ symmetrizers with condition numbers $10^{10}, 10^{11}, 10^{14}, 10^{15}, 10^{18},$ and $10^{19}$,  ranks below 23  (four times) and  full rank 23 only twice.\\[1.5mm]
The above examples all deal with sparse matrices $A_{n,n}$. We have also transformed each of these matrices separately  via  random dense unitary matrix similarities to dense form and then applied AZNN using the same parameter settings to obtain quite similar results numerically.\\[1mm]
 The main cost of AZNN methods for matrix square root and matrix symmetrizer computations lies in the necessary transformation of the respective $n$ by $n$ matrix model equations to Kronecker $n^2$ by $n^2$ matrix form and solving 
 the ensuing $n^2$ by $n^2$ linear system for each $t_k$. Note that $35^2 = 1225$ and the cost of solving each Kronecker matrix linear system now is $O(1225^3/3) = O(1.8383\cdot 10^9 / 3) \approx O(6 \cdot 10^8)$. Here the Kronecker system matrices of the unknown $X(t_k)$ both have the form $(X^T(t_k) \otimes I_n) + (I_n \otimes X(t_k))$, see \cite[Sections 2 (VIII) and  (VII) start-up)]{FUsurveyZNN}, whose sparsity structure may be exploited for sparse $A$ and $X$ and thus help us to expedite the linear system solves significantly. We have not  studied this idea any further.\\[2mm]
 Our last static matrix symmetrizer problem from \cite[p. 617 formula (9) through p. 619]{DU16} deals with the singular 2 by 2 matrix for varying $0 \leq \alpha \leq 1$\\[-5mm]
 $$ A = \bp 0& 1\\0 & \al\ep \ .
 $$
\hspace*{-.7mm}\includegraphics[height=69mm]{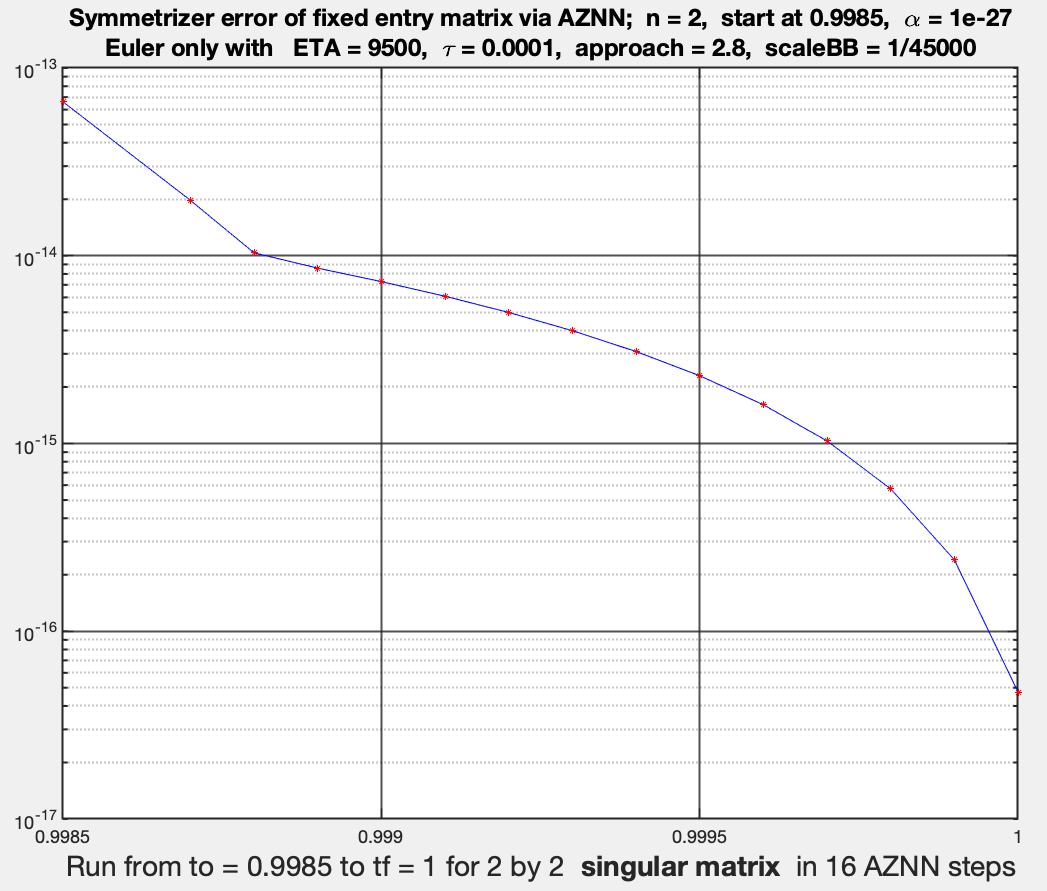} \vspace*{-2mm}\includegraphics[height=69mm]{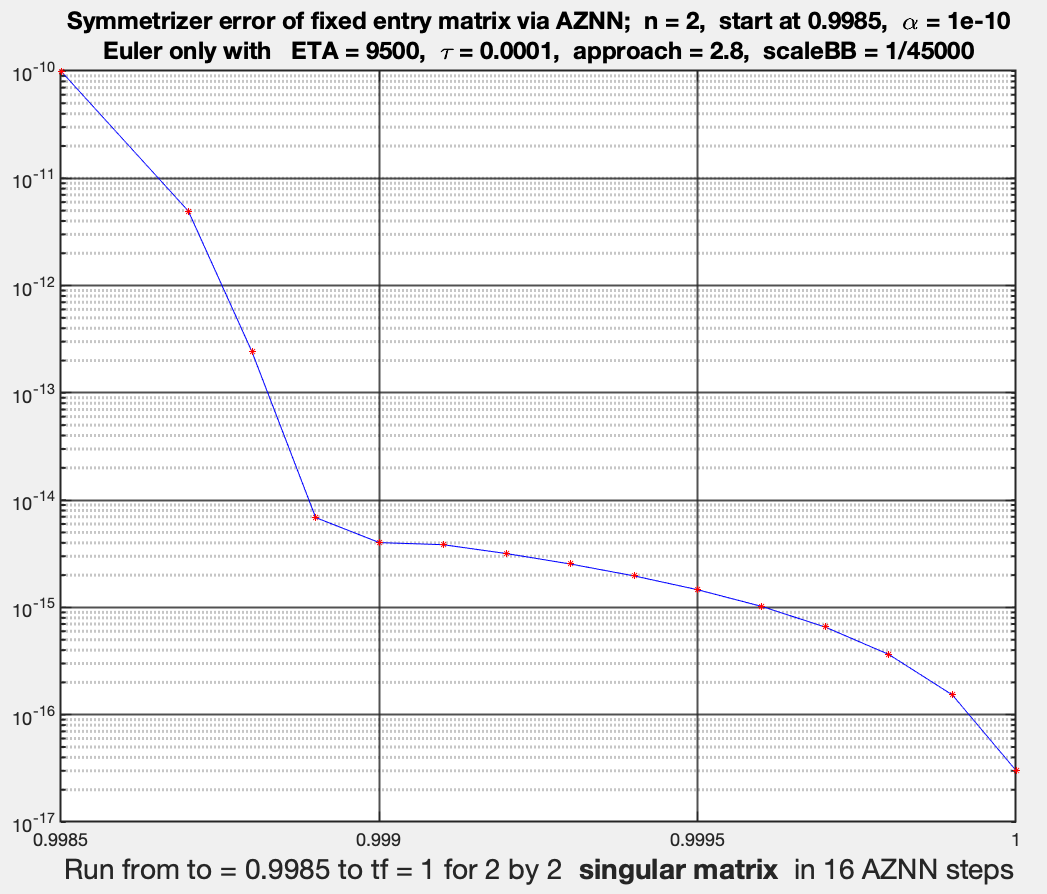}
\\[3.5mm]
\hspace*{-0.1mm}\includegraphics[height=69.5mm]{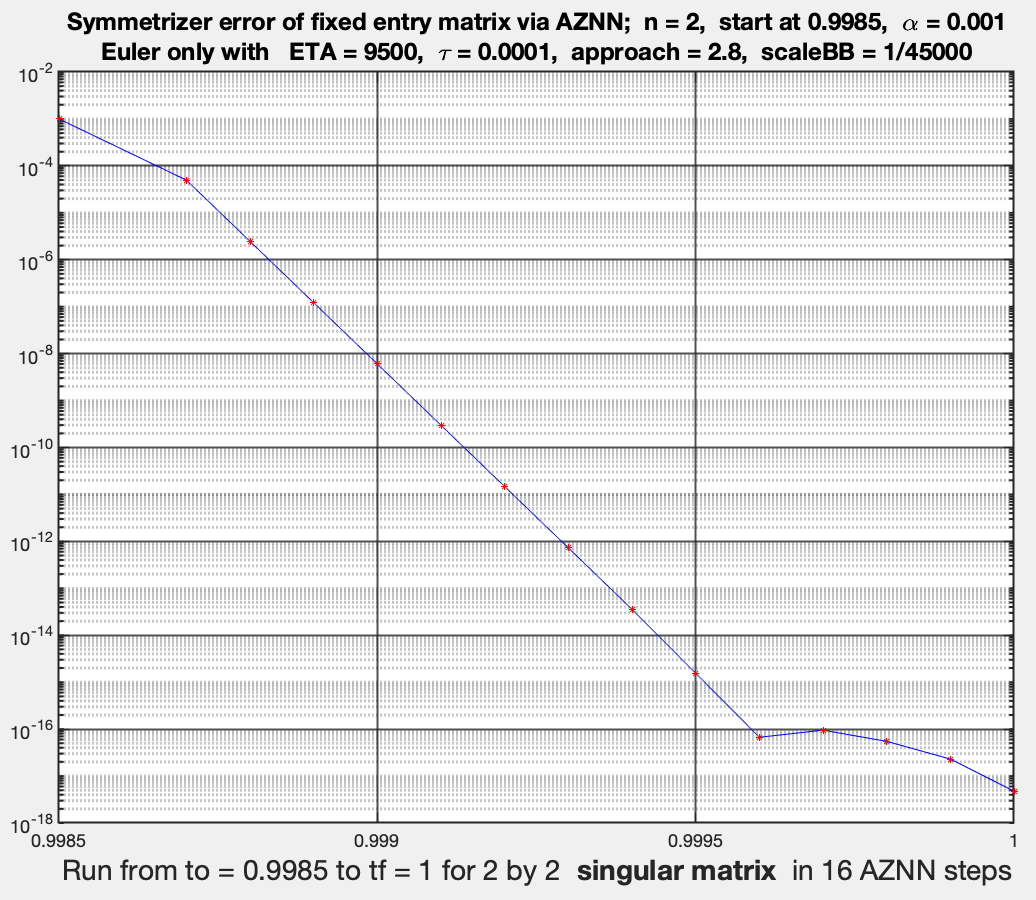} \hspace*{0.5mm}\includegraphics[height=69.5mm]{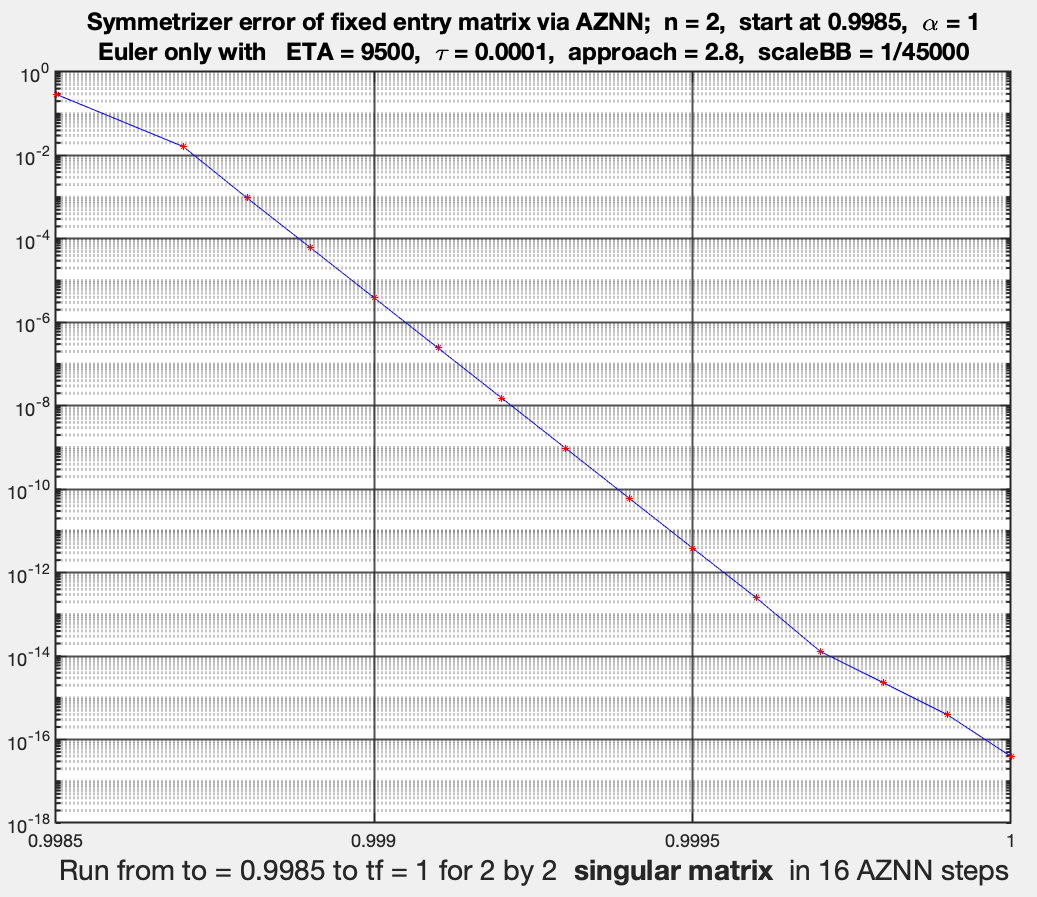}
\\[1mm]
\hspace*{73mm} Figure  12\\[1.5mm]
Here we receive excellent full rank results with symmetrizer condition numbers below 2.7 and relative symmetrizer errors well below the machine constant via AZNN for the identical parameter settings as used before.\\[1mm]
Our static matrix symmetrizer method is  AZNNsingleMatrSymmEuler3.m. It uses Matlab's 'gallery' matrices gallery('frank',35) or gallery('kahan',35) and Jordan block matrices of varying dimensions built in  AJordunisym.m of \cite{FUMatlabAZNN}.   \\[-8mm] 

\section*{Conclusions}

\vspace*{-0mm}
We have shown how to obtain faster more accurate  solutions for time-varying matrix problem through an adapted version AZNN of basic Zhang Neural Network (ZNN), where we define the exponential decay constants $ \eta$ in the start-up phase and in the ZNN iterations phase independently as the problem requires and also free the start-up phase to last as long as is needed to ensure a smooth transition between AZNN's two phases.\\[1mm]
 In the process we have also learnt how and when to use 'astronomically high' $\eta$ settings for improved theoretical results. These results were obtained when investigating  two complicated time-varying $n$ by $n$ matrix factorizations that require a model equation switch to  Kronecker matrix form due to their very formulation.\\[1mm]
For $n$ by $n$ matrix flows $A(t)_{n,n}$ and the unknown matrix $X(t)$,  Kronecker based matrix problems increase the linear equations solving dimensions for ZNN iterates from $O(n^3/3)$ to $O((n^2)^3/3) = O(n^6/3)$ operations. This slows down time-varying matrix factorizations from typically fractions of a second  to possibly dozens of seconds. But a study of the sparsity structure of the Kronecker matrices $(X^T(t) \otimes I_n) + (I_2 \otimes X(t))$ may improve the speed of time-varying matrix flow factorizations in the future. And really high dimensions $n$ are not generally encountered even in complex biological and chemical reaction models as demonstrated in \cite{EUA07} in numerous applications.\\[1mm]
This paper is not about {\em Theory} or the numerical analysis of time-varying matrix computations, which -- at this moment -- has not even begone for Zhang Neural Networks. Here we rather have shown some phenomena  and behavior that  appear with AZNN methods when altering the set of parameters of ZNN in each of its phases experimentally. One commonly noted event with basic ZNN  
is the relative constancy of the product $h = \eta \cdot \tau $ for ZNN over wide ranges of $\tau$ and $\eta$ for any one problem. Here is a list of $\eta $ and $\tau$ data for all our experiments which -- for AZNN -- puts this assumption to rest with wide $h$ variations and factors of 2 to 10 (Figure 6 left) and even 110 (Figure 3) for identical problems.\\[-4mm]
\begin{center}
\begin{tabular}{l||cc||cc}
&\multicolumn{2}{c||}{Left side graph}&\multicolumn{2}{c}{Right side graph}\\[0.2mm] \hline 
Figure 1 & \multicolumn{2}{c||}{h = 0.027}&\multicolumn{2}{c}{h = 0.056}\\[0.2mm]
(t.-v. square root)&&&&\\ \hline
Figure 3&\multicolumn{2}{c||}{h = 3.2}&\multicolumn{2}{c}{h = 0.029}\\
(t.-v. square root)&&&&\\ \hline \hline
Figures 4, 5 & \multicolumn{4}{c}{start \ \ \ \ \ middle \ \ \ \ \  end}\\
(t.-v. symmetrizer)& \multicolumn{4}{c}{ h = 28.4 \ \ \   h = 64 \ \ \ \ h = 56}\\ \hline
Figures 4, 5 & \multicolumn{4}{c}{start \ \ \ \ \ middle \ \ \ \ \  end}\\
(t.-v. symmetrizer)& \multicolumn{4}{c}{ h = 28.4 \ \ \   h = 64 \ \ \ \ h = 56}\\ \hline
Figure 6& start & end& start & end \\
(t.-v. symmetrizer)&h = 0.162 &h = 0.016&h = 0.162 &h = 0.084\\ \hline
Figure 7& start & end& start & end \\
(t.-v. symmetrizer)&h = 0.162&h = 0.018&h = 0.162 &h = 0.084\\ \hline
Figure 8&\multicolumn{4}{c}{start \ \ \ \ \ \ \ \ \ \ \ \ end}  \\
(t.-v. symmetrizer)&\multicolumn{4}{c}{h = 0.162 \ \ \ \ \ h = 0.018}\\ \hline
\hline
Figure 9&\multicolumn{2}{c||}{(Euler)}&\multicolumn{2}{c}{(Euler)}  \\
(static symmetrizer)&\multicolumn{2}{c||}{h = 0.29} &\multicolumn{2}{c}{h = 0.37}\\ \hline
Figure 10&\multicolumn{2}{c||}{(Euler)}&\multicolumn{2}{c}{(Euler)}  \\
(static symmetrizer)&\multicolumn{2}{c||}{h = 0.79} &\multicolumn{2}{c}{h = 0.67}\\ \hline
Figure 11&\multicolumn{4}{c}{(Euler)} \\
(static symmetrizer)&\multicolumn{4}{c}{h = 0.95} \\ \hline

Figure 12&\multicolumn{2}{c||}{top and bottom left graphs}&\multicolumn{2}{c}{top and bottom right graphs}  \\
(static symmetrizer)&\multicolumn{2}{c||}{(Euler)}&\multicolumn{2}{c}{(Euler)}  \\
&\multicolumn{2}{c||}{with $\alpha$ at $10^{-27}$ and $10^{-3}$}  
&\multicolumn{2}{c}{with $\alpha$ at $10^{-10}$ and  1} \\ 
&\multicolumn{2}{c||}{h = 0.95} &\multicolumn{2}{c}{h = 0.95}\\ \hline \hline
\end{tabular}\\[2mm]
Table of $h = \eta \cdot \tau$ values for Figures 1 through 12\\[2mm]
 Table 2
\end{center}

\vspace*{0mm}

\noindent
\centerline{{[} .. /box/local/latex/AZNN22/AZNN.tex] \quad \today }\\[1mm]

\noindent
20 image files in 12 figures :\\[1mm]

AZNNeta145.png\\
AZNNeta26.png\\
AZNNtau26large.png\\
AZNNsqrootsextended1houra.png\\
AZNNtvmatrSymm51houra.png\\

AZNNtvmatrSymm58hoursa.png\\
AZNNwiggle508.png\\
AZNNwiggle5422.png\\
AZNNnonwiggle542.png\\
AZNNnowiggle5098hours.png\\

\vspace*{-49.2mm}
\hspace*{60mm}
\begin{minipage}{55mm}{AZNNKahan2035.png\\
AZNNKahanlgpar35.png\\
AZNNFrank35.png\\
AZNNnonwiggle509.png\\
AZNNCutlgpar23.png\\

AZNN2by2alm27.png\\
AZNN2by2alm10.png\\
AZNN2by2alm3.png\\
AZNN2by2al1.png\\
AZNNFranklgpar35.png}
\end{minipage}\\[3mm]

Two  tables : \\[1mm]

Table1Fig4AZNN.PNG\\
Table 2

\end{document}